\documentclass[11pt, reqno]{amsart}

\usepackage[margin=1.05in]{geometry}
\usepackage{amsmath,amssymb,amsthm}
\usepackage{parskip}
\usepackage{enumerate}
\usepackage{tikz-cd}
\usepackage{graphicx}
\usepackage{caption}
\usepackage{subcaption}
\usepackage{tikz}
\usepackage{array}

\usepackage[colorinlistoftodos, disable]{todonotes}

\todostyle{Ric}{color=cyan}
\todostyle{Mig}{color=orange}
\todostyle{Mar}{color=yellow}

\newtheorem{theorem}{Theorem}[section]
\newtheorem{proposition}[theorem]{Proposition}
\newtheorem{lemma}[theorem]{Lemma}
\newtheorem{corollary}[theorem]{Corollary}
\newtheorem*{mainthm}{Main Theorem}
\theoremstyle{remark}
\newtheorem{remark}[theorem]{Remark}

\newcommand{\OO}{\mathbb O}
\newcommand{\CC}{\mathbb C}
\newcommand{\HH}{\mathbb H}
\newcommand{\RR}{\mathbb R}
\newcommand{\CaP}{\mathbb{O}\mathbb{P}^2}
\newcommand{\HP}{\mathbb{HP}}
\newcommand{\RP}{\mathbb{RP}}
\newcommand{\CP}{\mathbb{CP}}

\newcommand{\sph}{\mathbb{S}}

\newcommand{\tr}{\operatorname{tr}}
\newcommand{\id}{\operatorname{id}}
\newcommand{\spec}{\operatorname{spec}}

\newcommand{\fol}{\mathcal{F}}

\newcommand{\Sp}{\mathsf{Sp}}
\newcommand{\Spin}{\mathsf{Spin}}
\newcommand{\SU}{\mathsf{SU}}
\newcommand{\SO}{\mathsf{SO}}
\newcommand{\U}{\mathsf{U}}
\newcommand{\GG}{\mathsf{G}}
\newcommand{\Atype}[1]{\mathsf{A}_{#1}}
\newcommand{\Ctype}[1]{\mathsf{C}_{#1}}

\tikzset{dn/.style={circle,draw,inner sep=2.2pt}}
\newcommand{\nd}[2]{\node[dn,label={[font=\scriptsize,inner sep=1pt]above:#2}] (n#1) at (#1,0) {};}
\newcommand{\sedge}[2]{\draw (n#1)--(n#2);}
\newcommand{\dedge}[2]{%
  \draw ([yshift=1.1pt]n#1.east)--([yshift=1.1pt]n#2.west);
  \draw ([yshift=-1.1pt]n#1.east)--([yshift=-1.1pt]n#2.west);}
\newcommand{\dyn}[1]{\begin{tikzpicture}[x=1.15cm,baseline=(n0.base)]#1\end{tikzpicture}}

\author[M. Dom\'inguez-V\'azquez]{Miguel Dom\'inguez-V\'azquez}
\address{Universidade de Santiago de Compostela $-$ CITMAga\newline
\indent Department of Mathematics\newline
\indent Rúa Lope Gómez de Marzoa, s/n, \newline
\indent 15782 Santiago de Compostela, Spain}
\email{miguel.dominguez@usc.es}

\title{Polar foliations on the Cayley plane}
\author[R. A. E. Mendes]{Ricardo A. E. Mendes}
\address{University of Oklahoma\newline
\indent Department of Mathematics\newline
\indent 601 Elm Ave\newline
\indent Norman, OK, 73019-3103, USA}
\email{ricardo.mendes@ou.edu}

\author[M. Radeschi]{Marco Radeschi}
\address{Universit\`a degli Studi di Torino\newline
\indent Departimento di Matematica ``G. Peano''\newline
\indent Via Carlo Alberto, 10\newline
\indent 10123 Torino (TO), Italy}
\email{marco.radeschi@unito.it}

\thanks{The first author has been supported by grant PID2025-175773NB-I00 funded by MICIU/AEI/10.13039/501100011033. The second author has been supported by the NSF grant DMS-2506409. }

\subjclass[2020]{53C12, 53C35}
\keywords{Polar foliation, isoparametric foliation, Cayley projective plane}

\begin{document}\maketitle
\todo[Ric]{changed the order of the authors so it's alphabetical. Changed grant acknowledgements accordingly, and added mine.}
\begin{abstract}\noindent
We classify polar foliations on the Cayley projective plane, completing the classification of polar foliations on compact rank one symmetric spaces.
\end{abstract}

\section{Introduction}

A polar foliation $\fol$ on a Riemannian manifold $M$ is a singular Riemannian foliation on $M$ such that for any $p\in M$ there is a totally geodesic submanifold $\Sigma$ of $M$ containing $p$ and intersecting all leaves of $\fol$ perpendicularly. A singular Riemannian foliation on $M$ is called homogeneous if it is the family of orbits of a Lie group isometric action on $M$. 

The main result of this paper is the classification of (not necessarily homogeneous) polar foliations on the Cayley projective plane $\CaP$ up to congruence (i.e.\ up to global isometry of $\CaP$), thus completing the classification of polar foliations on the compact rank one symmetric spaces (CROSSes).

\begin{mainthm}\todo[Ric]{Removed the number ``1''. There is only one main theorem.}
	Let $\fol$ be a polar foliation on the Cayley projective plane $\CaP$. Then, up to congruence, $\fol$ is one of the foliations in Table \ref{T:cases} below. In particular, $\fol$ is either a homogeneous foliation, or congruent to a certain inhomogeneous foliation $\fol_c$ of codimension $2$.
\end{mainthm}

\renewcommand{\arraystretch}{2.1}
\begin{table}[ht]
\centering
\begin{tabular}{|c|c|c|c|c|c|}
\hline
Case&$\mathcal{F}$ & $\operatorname{codim}$ & Lifted Weyl group & Marked diagram & Leaf space \\
\hline\hline
$1a$&$\Spin(9)$ & 1 & $\tilde{\mathsf{A}}_1$ &
  \dyn{\nd{0}{15}\nd{1}{7}\sedge{0}{1}
       \node[font=\scriptsize] at (0.5,0.25) {$\infty$};}
  & $[0,\tfrac{\pi}{2}]$ \\
\hline
$1b$&\begin{tabular}{c}
$\Sp(3)\cdot\Sp(1)$\\ $\Sp(3)\cdot\U(1)$\\ $\Sp(3)$
\end{tabular} & 1 & $\tilde{\mathsf{A}}_1$ &
  \dyn{\nd{0}{7}\nd{1}{4}\sedge{0}{1}
       \node[font=\scriptsize] at (0.5,0.25) {$\infty$};}
  & $[0,\tfrac{\pi}{4}]$ \\
\hline\hline
$2$&$\Spin(8)$ & 2 & $\Atype{1}\times\Atype{1}\times\Atype{1}$ &
  \dyn{\nd{0}{7}\nd{1}{7}\nd{2}{7}}
  & $\Delta(2,2,2)$ \\
\hline\hline
$3a$&$\Spin(7)\cdot\SO(2)$ & 2 & $\Atype{1}\times\Ctype{2}$ &
  \dyn{\nd{0}{7}\nd{1}{6}\nd{2}{1}\dedge{1}{2}}
  & $\Delta(2,2,4)$ \\
\hline
$3b$&$\SU(4)\cdot\SU(2)$ & 2 & $\Atype{1}\times\Ctype{2}$ &
  \dyn{\nd{0}{7}\nd{1}{5}\nd{2}{2}\dedge{1}{2}}
  & $\Delta(2,2,4)$ \\
\hline
$3c$&$\mathcal{F}_c$ inhomog. & 2 & $\Atype{1}\times\Ctype{2}$ &
  \dyn{\nd{0}{7}\nd{1}{4}\nd{2}{3}\dedge{1}{2}}
  & $\Delta(2,2,4)$ \\
\hline\hline
$4a$&$\SU(3)\cdot\SU(3)$ & 2 & $\Ctype{3}$ &
  \dyn{\nd{0}{2}\nd{1}{2}\nd{2}{3}\sedge{0}{1}\dedge{1}{2}}
  & $\Delta(2,3,4)$ \\
\hline
$4b$&$\SO(3)\cdot\GG_2$ & 2 & $\Ctype{3}$ &
  \dyn{\nd{0}{1}\nd{1}{1}\nd{2}{5}\sedge{0}{1}\dedge{1}{2}}
  & $\Delta(2,3,4)$ \\
\hline
\end{tabular}
\caption{Polar foliations on $\CaP$. The lifted Weyl group $\tilde\Gamma$ of such a foliation is recalled in \S \ref{S:polar-fol}, and it is used to group up the cases. The marked diagram is the Dynkin diagram of $\tilde{\Gamma}$, with weights given by the multiplicities of the codimension 1 strata fixed by the corresponding reflections (the notion of multiplicity is recalled in \S \ref{S:SRFs}). $\Delta(p,q,r)$ denotes a spherical triangle of curvature 1, with angles $(\tfrac{\pi}{p},\tfrac{\pi}{q},\tfrac{\pi}{r})$.}
\label{T:cases}
\end{table}

The unique inhomogeneous polar foliation $\fol_c$, and its properties, are described in Proposition~\ref{P:nonhom-fol}, whereas the homogeneous ones are already classified~\cite{PT,GK16}.
%

The history of the study of polar foliations is long and fascinating. By a recent result by Liu and Radeschi~\cite{LR22}, polar foliations on compact symmetric spaces are \emph{isoparametric}, i.e.\ the mean curvature of leaves is a basic vector field (or equivalently, its flow moves leaves to leaves). Codimension one isoparametric foliations were first studied in space forms by Segre and Cartan in the 30s. It turned out that the case of spheres was the hardest: initiated by Cartan, the classification was completed after the concerted efforts of multiple mathematicians spanning a century: Cartan \cite{Car38, Car39}, M\"unzner \cite{Mun80,Mun81}, Ferus-Karcher-M\"unzner \cite{FKM}, Abresch \cite{Ab83}, Dorf\-meister-Neher \cite{DN85}, Stolz \cite{Sto99},  Cecil-Chi-Jensen~\cite{CCJ}, Immervoll \cite{Imm08}, Miyaoka \cite{Mi13}\todo[Mig]{I cited Abresch paper (instead of thesis), and Miyaoka's Annals paper. Even if probably incorrect, most of its results are correct. We may cite Siffert as well}, Chi \cite{Chi20}, Qian-Tang \cite{QT26}, and many others. We refer to \cite{Chi20, QT26, Th:sp} for further information and bibliography.
 
Isoparametric foliations of arbitrary codimension on spaces forms were introduced by Terng~\cite{Te:jdg} in the 80s. The main case is again that of spheres, where in higher codimension the picture is rigid. Indeed, Thorbergsson~\cite{Th:buildings} showed that irreducible \todo[Mig]{Full deleted: I think needed for submanifolds, but not for foliations, right?}
polar foliations of codimension at least two in round spheres are homogeneous. An alternative proof was obtained by Olmos~\cite{Olmos:jdg}.

Later, Terng and Thorbergsson~\cite{TT:jdg} extended the theory to symmetric spaces, and Heintze-Liu-Olmos~\cite{HLO} to general Riemannian manifolds. Years later, Lytchak \cite{Lyt14} proved that polar foliations of codimension at least three on irreducible compact symmetric spaces are hyperpolar (that is, polar with flat sections) unless the space has rank one. Irreducible hyperpolar foliations were known to be homogeneous, by a result by Christ~\cite{christ}, and homogeneous polar foliations on compact symmetric spaces of rank at least two
have been classified by Kollross and Lytchak~\cite{Kollross:jdg, KL:blms}. Very recently, Wiesendorf~\cite{Wie26} extended Lytchak's hyperpolarity result to codimension two polar foliations.

The difficulty is thus concentrated in compact rank one symmetric spaces (CROSSes). Polar actions on complex, quaternionic and octonionic projective spaces were classified by Podestà and Thorbergsson \cite{PT} (with an omission detected by Gorodski and Kollross~\cite{GK16}). The classification of (not necessarily homogeneous) polar foliations on complex and quaternionic projective spaces $\mathbb{CP}^n$ and $\mathbb{HP}^n$ was accomplished in a series of papers  by Domínguez-Vázquez, Gorodski, and Kollross~\cite{Do:tams, DG:tohoku, DK:mathz}. 
 

The analysis of polar foliations on the Cayley plane $\CaP$ is fundamentally different from the case of the other projective spaces $\mathbb{KP}^n$, $\mathbb{K}\in\{\CC, \HH\}$. In these cases, the existence of a Hopf fibration $\sph^m\to \mathbb{KP}^n$ suggests lifting the foliation to the sphere, understanding which foliations can occur, and how they project back to $\mathbb{KP}^n$ (congruent foliations on $\sph^m$ can produce incongruent foliations on $\mathbb{KP}^n$). 
The absence of a Hopf fibration over the Cayley plane $\CaP$ requires a completely different strategy. Ours is based on the following steps.  First, by analysing the possible geometry of the quotient, one gets to a natural division into different cases. For each case, the main step consists in proving that there is always a totally geodesic singular leaf $L_*$, and in identifying the specific model. By Wolf's classification of totally geodesic submanifolds \cite{Wolf}, this fixes $L_*$ up to congruence. The rest of the congruence is established by analysing the infinitesimal foliation around $L_*$. In order to prove the existence of a totally geodesic leaf, we analyse the evolution of certain tensors, including the shape operator, along horizontal geodesics, and connecting their properties to the focal sets of leaves.

After completing the present work, the authors were informed of the recent work by Wiesendorf \cite{Wie26} which, apart from the hyperpolarity result for codimension two polar foliations mentioned above, also recovers the homogeneity of polar foliations in item cases $4a$, $4b$ of Table~\ref{T:cases}. \todo[Ric]{Separated this paragraph from the next.}

Combining our classification with the results mentioned above
it follows that\todo[Mig]{Let me try to explain better my doubts about including this claim. One also has the open question of the existence of inhomogeneous reducible foliations on higher rank spaces. Since for $\CaP$ this is not trivial, that is not expected to be easy, either. And in order to just state that the combination of all known results reduces de classification of irreducible foliations, then we don't need our result (Wiesendorf's is enough).}
a full classification of polar foliations on simply connected  symmetric spaces of compact type reduces to a classification of codimension one (polar) foliations on spaces of rank higher than one.

The paper is structured as follows. In Section~\ref{sec:preliminaries} we recall some basic facts about $\CaP$ and some concepts and tools from the theory of singular Riemannian foliations. In Section~\ref{subsec:lagrangian} we define the tensors along horizontal geodesics that we will need, and their relation to focal distances. In Section~\ref{S:codim-1} we classify polar foliations $(\CaP,\fol)$ of codimension 1, and recover cases $1a$, $1b$ of Table \ref{T:cases}. In Section~\ref{S:codim-2} we classify polar foliations of codimension 2, and recover the remaining cases in Table~\ref{T:cases}: case $2$ is obtained in Proposition~\ref{prop:Spin8}, cases $3a$, $3b$, $3c$ in Theorem~\ref{T:pi/4}, and cases $4a$, $4b$ in Theorem~\ref{T:irred}. 

\textbf{AI disclosure.}
The research strategy and fundamental ideas were set by the authors. The authors used Claude Opus 5.5 to explore possible approaches to the technical analysis of the codimension one case, and in producing the estimates that were crucial in the classification of cases of type $4$ in Table~\ref{T:cases}. All results and proofs were written and verified by the authors. The manuscript contains no AI-generated text.

\subsection*{Acknowledgements}
The authors would like to thank A. Kollross for insightful discussions, on a previous approach to the Main Theorem.

\section{Preliminaries}\label{sec:preliminaries}
\subsection{The Riemannian geometry of the Cayley plane}\label{S:CaP-geometry}
Let $\CaP=\mathsf{F}_4/\Spin(9)$ carry its symmetric metric with sectional curvature in $[1,4]$. Then, its geodesics are closed of length $\pi$. We denote by $R$ its curvature tensor, with the usual convention that $\langle R(y,x)x,y\rangle$ is the sectional curvature of the plane spanned by two orthonormal tangent vectors $x,y$.

Let $p\in \CaP$ and $x\in T_p\CaP$ be a unit tangent vector. The Jacobi operator $R_x=R(\cdot,x)x$  is a self-adjoint endomorphism of $x^\perp$ with eigenvalues $1$ and $4$. Denote by $E_1(R_x)$
and $E_4(R_x)$ the corresponding eigenspaces, which satisfy $\dim E_1(R_x)=8$ and $\dim E_4(R_x)=7$. Let $\ell(x)$ denote the octonionic line (or Cayley line) through $p$ tangent to $x$, defined as the tangent space to a totally geodesic $\mathbb{OP}^1=\mathbb{S}^8(4)$ tangent to $x$. Then
\[
\ell(x)=\mathbb{R}x\oplus E_4(R_x), \qquad \ell(x)^\perp=E_1(R_x).
\]
From now on, by $P_V$ we denote the orthogonal projection map onto a subspace $V$ of the tangent space to $\CaP$ at some point.

The following result is elementary and well known: we report it here for easy reference later on.

\begin{lemma}\label{lem:cayley}
Let $x,y\in T_p\CaP$ be orthonormal. Then
\begin{enumerate}[\rm(1)]
\item $\langle R_xy,y\rangle=1+3\,|P_{\ell(x)}y|^2$. In particular $\langle R_xy,y\rangle=1$ if and only if $y\perp\ell(x)$, and
      $\langle R_xy,y\rangle=4$ if and only if $y\in\ell(x)$.
\item $y\in\ell(x)$ if and only if $\ell(y)=\ell(x)$, whereas $y\perp\ell(x)$ if and only if $\ell(y)=\ell(x)^{\perp}$.
\item A subspace of $T_p\CaP$ all of whose $2$-planes have curvature $1$ has dimension $\le2$.
\end{enumerate}
\end{lemma}

\begin{proof}
(1) On $x^\perp$ we have $R_x=P_{E_1(R_x)}+4P_{E_4(R_x)}$, so
$\langle R_xy,y\rangle=|P_{E_1(R_x)}y|^2+4|P_{E_4(R_x)}y|^2=1+3|P_{E_4(R_x)}y|^2$, and
$P_{E_4(R_x)}y=P_{\ell(x)}y$ because $y\perp x$.

(2) The first equivalence holds because a Cayley line through $p$ is determined by any of its tangent
directions. For the second, $\Spin(9)$ is transitive on the sphere of unit vectors, so we may use the model
$T_q\CaP=\mathbb{O}\oplus\mathbb{O}$ with $x=(1,0)$ and $\ell(x)=\mathbb{O}\oplus0$; then $y\perp\ell(x)$
means $y\in0\oplus\mathbb{O}$, which is again an octonionic line.

(3) If $x,y,z$ are orthonormal with all sectional curvatures equal to $1$, then $y,z\perp\ell(x)$, so
$\ell(y)=\ell(x)^\perp\ni z$ and $\langle R_yz,z\rangle=4$: a contradiction.
\end{proof}

\subsection{Singular Riemannian foliations}\label{S:SRFs}

We summarize some facts about singular Riemannian foliations, and refer to~\cite{Rad26} for a detailed exposition.
A \emph{singular Riemannian foliation} is a partition $(M,\fol)$ of a Riemannian manifold $M$ into connected, smooth, injectively immersed submanifolds, called leaves, which are locally equidistant to each other and whose tangent spaces are spanned by a family of smooth vector fields. The local equidistance implies that any geodesic starting perpendicular to a leaf stays perpendicular to all leaves it meets. We will refer to them as \emph{horizontal geodesics}.

Given a point $p\in M$, we denote by $L_p$ the leaf through $p$. The dimension $\dim L_p$ might depend on $p$: the maximum dimension will be called $\dim \fol$, the \emph{dimension of $\fol$}. Similarly, $\dim M-\dim \fol$ is called the \emph{codimension~of~$\fol$}.

A point $p$ is called \emph{regular} if $\dim L_p=\dim \fol$, and \emph{singular} otherwise. For every $d=0,\ldots, \dim \fol$, we define the \emph{$d$-dimensional strata of $\fol$} as the connected components $M_d^i$ of the set
\[
M_d=\{p\in M : \dim L_p=d\}.
\]
Every stratum is a smooth (possibly noncomplete) submanifold, and for any stratum $M_d^i$, one has that $\overline{M_d^i}\setminus M_d^i$ is a union of strata of strictly lower dimensions. The set $M_{\dim \fol}$ is connected, open, dense, and of full measure in $M$. We call this stratum the \emph{regular stratum of $\fol$} and denote it by $M^{reg}$. For any other stratum $M_d^i$, we call the difference $\dim \fol- d$ the \emph{multiplicity of $M_d^i$}. Similarly, given a leaf $L$, we define the \emph{multiplicity of $L$} as the difference $\dim \fol-\dim L$.

\subsection{Polar foliations}\label{S:polar-fol}
Most of the content of this subsection can be found in \cite{AT13, Ale11, Toben} and references therein. A singular Riemannian foliation $(M,\fol)$ is called \emph{polar} if through each point there exists a submanifold $\Sigma\subset M$ that meets all leaves of $\fol$, and it does so perpendicularly. Such a submanifold is called a \emph{section} and it is necessarily totally geodesic.

If $(M,\fol)$ is a polar foliation on a simply connected manifold, then the leaves of $\fol$ are closed \cite[Corollary 1.4]{Ale11}. In this case, the leaf space $X=M/\fol$ inherits the structure of a Hausdorff metric space, such that the projection map $\sigma\colon M\to X$ is a submetry (i.e. it takes metric $r$-balls to metric $r$-balls).

If $(M,\fol)$ is a polar foliation with $M$ compact and simply connected, then given a section $\Sigma$ of $\fol$ there exists a finite group ${\Gamma}$, called \emph{Weyl group of $\Sigma$}, acting on $\Sigma$ by isometries with the following properties:
\begin{enumerate}
\item $\Gamma$ is generated by \emph{reflections}, i.e.\ involutions $r_i\colon \Sigma\to \Sigma$ that fix a codimension-1 submanifold, called \emph{wall of $r_i$} \cite[Theorem 1.1]{Ale11}.
\item For any point $p\in \Sigma$ the orbit $\Gamma\cdot p$ coincides with $L_p\cap \Sigma$. In particular, the inclusion $\Sigma\to M$ induces an isometry
\[
\Sigma/\Gamma\to M/\fol.
\]
Therefore, in this case the quotient $M/\fol$ is isometric to a \emph{good orbifold}, i.e.\ a quotient of a smooth manifold by a discrete group acting properly on it.
\end{enumerate}

\todo[Mar]{added info about lifted Weyl group}
Furthermore, this Weyl group lifts to a reflection group $\tilde{\Gamma}$ on the universal cover $\tilde{\Sigma}$ of the section, which normalizes the deck transformation group of $\tilde{\Sigma}\to \Sigma$. We will call this the \emph{lifted Weyl group}.

In general, every stratum $M_d^i$ projects, via $\sigma\colon M\to M/\fol$, to a smooth orbifold. However, if $\fol$ is polar and $M$ is simply connected, the following structure holds:
\begin{enumerate}
\item The leaf space $X$ is isometric to (the closure of) a fundamental domain in $\Sigma$ for the $\Gamma$-action (these fundamental domains are called \emph{Weyl chambers}). This follows from the fact that $W$ acts on $\Sigma$ by reflections.
\item By \cite[Theorem 1.6]{Lyt10} the regular stratum $M^{reg}$ projects to the interior of $X$, while the singular strata project to the boundary of~$X$.
\end{enumerate}

{\bf In our case} we will be interested in polar foliations $(\CaP,\fol)$ on the Cayley plane. By \cite[Proposition 9.2]{Lyt14}, these either have codimension 1 or codimension 2.

If the codimension is 1, a section $\Sigma$ must be a closed geodesic $\sph^1$ (of length $\pi$ by the normalization of $\CaP$ from \S\ref{S:CaP-geometry}). Since $\CaP/\fol\simeq \sph^1/\Gamma$ must be simply connected (see e.g. \cite[Lemma 2.2]{Suk26}), it follows that $\Gamma\neq \{1\}$ is a reflection group on $\Sigma=S^1$ (hence a dihedral group $D_{2g}$) and $\CaP/\fol\simeq \Sigma/\Gamma$ is isometric to an interval $[0, \tfrac{\pi}{2g}]$. If the codimension is 2,  again \cite[Proposition 9.2]{Lyt14} proves that sections are congruent to a totally geodesic $\RP^2$ with constant curvature 1.

\subsection{The infinitesimal foliation at a point}\label{S:inf-fol}
The content of this subsection can be found in \cite{Mol86}, as well as in \cite{Rad26}. Given a singular Riemannian foliation $(M,\fol)$ and given $p\in M$, the normal space $\nu_p L_p$, with its Euclidean metric $g_p$, admits a singular Riemannian foliation $\nu_p\fol$ with the following properties:
\begin{enumerate}
\item The exponential map $\exp_p\colon \nu_pL_p\to M$ maps leaves of $\nu_p\fol$ to leaves of $\fol$ on a small ball around the origin.
\item The origin is a leaf. In particular, $(\nu_p L_p,\nu_p\fol)$ restricts to a singular Riemannian foliation $(\sph_p,\fol_p)$ on the unit sphere $\sph_p$ of $\nu_pL_p$.
\item If $(M,\fol)$ is polar with section $\Sigma$ through $p$, then $(\nu_pL_p, \nu_p\fol)$ is polar with section $T_p\Sigma$, and $(\sph_p,\fol_p)$ is polar with section $T_p^1\Sigma$.
\item Given $\sigma\colon M\to M/\fol$, then the (metric) tangent cone of $M/\fol$ at $\sigma(p)$ is isometric to a quotient of $\nu_pL_p/\nu_p\fol$ by a finite group of isometries acting on it. Furthermore, $\nu_pL_p/\nu_p\fol$ is isometric to the Euclidean cone over $\sph_p/\fol_p$.
\end{enumerate}

We will use these facts in the special case in which $M/\fol$ is a 2-dimensional polyhedron. Then, fixing a vertex $p_*\in M/\fol$ with internal angle $\theta$, the tangent cone $T_{p_*}M/\fol$ is the cone over a segment of length $\theta$. In particular, given a leaf $L$ projecting to $p_*$ and a point $p\in L$, the infinitesimal foliation $(\sph_p, \fol_p)$ is a codimension one foliation in a round sphere.

\subsection{Codimension one polar foliations on spheres}

Codimension one foliations in round spheres appear in the literature with the name of \emph{codimension one isoparametric foliations}, and the regular leaves take the name of \emph{isoparametric hypersurfaces}. Given a codimension one foliation $(\sph^n, \fol)$, every horizontal geodesic is a section. As in the case of the Cayley plane discussed above, since geodesics in $\sph^n$ are closed $\sph^1$'s of length $2\pi$, the Weyl group acting on one such section cannot be trivial, and it must be a dihedral group $D_{2g}$. Therefore, the leaf space $\sph^n/\fol$ is isometric to an interval of length $\pi/g$.

Letting $L_-$, $L_+$ be the singular leaves projecting to the endpoints of the interval, denote by $m_\pm=\dim\fol-\dim L_\pm$ their multiplicities. Two facts that we will use later on are:
\begin{itemize}
\item The fundamental result of M\"unzner \cite{Mun81}, that $g\in \{1,2,3,4,6\}$.
\item The identity proved in \cite[Satz C]{Mun80}, that
\begin{equation}\label{E:relation}
2(n-1)=g(m_++m_-).
\end{equation}
\end{itemize}

\section{Lagrangian families of $L$-Jacobi fields}\label{subsec:lagrangian}

Let $L$ be a submanifold of $\CaP$, fix $p\in L$ and let $W:=T_pL$ denote the tangent space and $\nu:=\nu_pL$ the normal space.

Fixing a unit vector $x\in \nu$, we abbreviate $H_x:=\nu\cap x^{\perp}$ and define
\[
\Pi_x^1:=P_{E_1(R_x)}=P_{\ell(x)^{\perp}},\qquad \Pi_x^4:=P_{E_4(R_x)}=P_{\ell(x)\cap x^{\perp}},
\]
and hence we have $\Pi_x^1+\Pi_x^4=\mathrm{id}_{x^{\perp}}$ and $R_x=\Pi_x^1+4\Pi_x^4$. 

For a unit $x\in\nu$, we consider the geodesic $\gamma_x(t)=\exp_p(tx)$ of $\CaP$. Since $\CaP$ is a symmetric space, both the Jacobi operator $R_{{\gamma}_x'}$ and its eigenspaces are parallel along $\gamma_x$. To simplify the formulas, we abuse the notation $\Pi_x^i$, and use it to denote $P_{E_i(R_{{\gamma}_x'})}$, for $i\in\{1,4\}$. With this in mind, we introduce the following operators on $({\gamma}_x')^{\perp}$:
\begin{equation}\label{E:OperatorsGivingJ}
F(t):=\cos t\,\Pi_x^1+\cos2t\,\Pi_x^4,\qquad
D(t):=\sin t\,\Pi_x^1+\tfrac{\sin2t}{2}\,\Pi_x^4 .
\end{equation}

Let $S_x$ be the shape operator of $L$ with respect to the unit normal $x$. Along the geodesic $\gamma_x$, the $L$-Jacobi fields, i.e.\ Jacobi fields arising from variations through horizontal geodesics starting from $L$, are characterized by the initial conditions (see~\cite[\S 10.2]{BCO})
\[
J_v(0)=w\qquad \textrm{and}\qquad J'_v(0)=-S_xw+y,
\]
where $v=(w,y)\in W\oplus H_x$ and the correspondence $v\mapsto J_v$ is linear. In particular,
the space $\Lambda$ of $L$-Jacobi fields (also called the Lagrangian family of $L$-Jacobi fields~\cite{Lyt09}) is a vector space isomorphic to $x^\perp$ via the map
\begin{align*}
x^\perp=W\oplus H_x&\longrightarrow \Lambda\\
v=(w,y)&\longmapsto J_v.
\end{align*}
For each $t\in\mathbb{R}$ we also denote by $J(t)$ the evaluation map $W\oplus H_x\longrightarrow T_{\gamma_x(t)}\CaP$ given by $v\mapsto J_v(t)$. 
Since $R_x$ is
parallel along $\gamma_x$, identifying vectors with their parallel translations along $\gamma_x$, one can easily check that
\begin{equation}\label{E:Jacfield}
J_v(t)=F(t)\,w-D(t)\,(S_xw-y) 
\end{equation}
satisfies the Jacobi equation $J_v''+R_{\gamma_x'}(J_v)=0$ and the required initial conditions.

Given a $t\in \RR$, the dimension of $\ker(v\mapsto J_v(t))$ counts the number of $L$-focal points at $t$, counted with multiplicity. We can compute this dimension in two different ways. From the ``foliation side'', assume $L$ is a leaf of a polar foliation $(\CaP,\fol)$ with leaf space $X=\CaP/\fol$ and canonical projection $\sigma\colon \CaP\to X$. Then by \cite[Section 3.2]{Lyt09}, the number of focal points of $L$ (counted with multiplicities) along a horizontal geodesic $\gamma_x$ passing through the interior of $X$, at time $t$, is the sum of:
\begin{enumerate}
\item A \emph{vertical contribution}, given by the difference $\dim L_0-\dim L_{\gamma_x(t)}$ between the dimension of a regular leaf $L_0$ of $\fol$ and the leaf through $\gamma_x(t)$.
\item A \emph{horizontal contribution}, given by the multiplicity of $\sigma(\gamma_x(t))$ as a conjugate point of $\sigma(\gamma_x(0))$, along the geodesic $\sigma\circ\gamma_x$.
\end{enumerate}

From the ``octonionic side'', the multiplicity has to do with the mutual position between $W=T_pL$ and the eigenspaces $E_1(R_x)$, $E_4(R_x)$ of the normal Jacobi operator. Thus, it will be important to consider the self-adjoint operators $B_x, C_x\colon W\to W$ given by
\[
B_x:=P_W\Pi_x^1\big|_W, \qquad\quad
C_x:=\id_W-B_x=P_W\Pi_x^4\big|_W,
\]
and define 
\[
Z_x:=\ker C_x=W\cap E_1(R_x),
\]
where the last equality holds because $\langle C_xw,w\rangle=|P_{E_4(R_x)}w|^{2}$.

\begin{remark}\label{R:C-and-curv}
If $w\in W$ and $x\in \nu$ are unit vectors, by Lemma \ref{lem:cayley} we have
\[
\langle C_xw,w\rangle=\tfrac13 (\langle R_xw,w\rangle-1)=\tfrac13 (\operatorname{sec}(x,w)-1).
\]
\end{remark}

\begin{proposition}\label{prop:pencil}
Let $L$ be a submanifold of $\CaP$, $p\in L$, $x\in \nu$ unit, and let $\Lambda$ denote the Lagrangian family of $L$-Jacobi fields along $\gamma_x$, with corresponding map $J(t)\colon x^\perp =W\oplus H_x\to T_{\gamma_x(t)}\CaP$.
\begin{enumerate}[\rm(a)]
\item For $t\in(0,\pi)$, $t\ne\pi/2$, set $c=\cot t$, and let $w\in W$. Then $J(t)(w,y)=0$ for some $y\in H_x$ if and only if $\big(c^{2}\id_W-c\,S_x-C_x\big)w=0$. In this case, $y$ is uniquely determined by $w$. In particular,
\begin{equation}\label{Equation-not-piover2}
\dim \ker J(t)=\dim\ker\big(c^2\id_W-cS_x-C_x\big).
\end{equation}
\item For $t=\pi/2$, $J_{(w,y)}(\pi/2)=0$ if and only if $w\in Z_x$ and $P_{Z_x}S_xw=0$. In this case, $y$ is determined by $w$
modulo $H_x\cap E_4(R_x)$. In particular:
\begin{align}\label{Gen-Equationpiover2}
\dim\ker J(\pi/2)=&\dim\bigl(H_x\cap E_4(R_x)\bigr)+\dim\ker\bigl(P_{Z_x}S_x|_{Z_x}\bigr).
\end{align}
\end{enumerate}
\end{proposition}

\begin{proof}
(a) By \eqref{E:Jacfield}, $J_{(w,y)}(t)=0$ if and only if $F(t)w=D(t)(S_xw-y)$. Decomposing this equation along $E_1\oplus E_4$ for $t\ne\pi/2,\pi$ gives:
\begin{align*}
(S_xw-y)_1=&\cot t\,w_1= c w_1,\\
(S_xw-y)_4=&\frac{2\cos2t}{\sin 2t}\,w_4=(\cot t-\tan t)\,w_4=(c-c^{-1})w_4,
\end{align*}
where we used the short-hand notation $w_1$, $w_4$ to denote $\Pi^1_xw$, $\Pi^4_xw$, respectively.
Therefore if $w$ is such that such a $y$ exists, then
\[
y=y_1+ y_4=\Pi^1_x(S_xw-cw)+\Pi^4_x\left(S_xw-(c-c^{-1})w\right).
\]
Hence, for $w\in W$ there exists a $y$ such that $(w,y)\in \ker J(t)$ if and only if
\begin{align*}
0=P_Wy=&P_W\Pi^1_x(S_xw-cw)+P_W\Pi^4_x\left(S_xw-(c-c^{-1})w\right)\\
=&B_x(S_xw-cw)+C_x\left(S_xw-(c-c^{-1})w\right)\\
^{(B_x+C_x=\id_W)}=&S_xw-cw+c^{-1}C_xw,
\end{align*}
and in this case $y$ is uniquely determined by $w$. Therefore, $J(t)(w,y)=0$ if and only if $w\in \ker(S_x-c\id_W+c^{-1}C_x)$. 
Multiplying by $-c=-\cot t$ gives the claim.

(b) At $t=\pi/2$, $F(\pi/2)=-\Pi_x^4$ and $D(\pi/2)=\Pi_x^1$, so $J_{(w,y)}(\pi/2)=0$ becomes
\[
w_4=0,\qquad(S_xw-y)_1=0.
\]
The first says $w\in W\cap E_1(R_x)=Z_x$.  The second says
$S_xw-y\in E_4(R_x)$, that is
\[S_xw\in H_x+E_4(R_x)=(W\cap E_1(R_x))^\perp\cap x^\perp=(\mathbb{R}x\oplus Z_x)^\perp.
\]
Since $S_xw$ lies in
$W$, the condition is $S_xw\perp Z_x$, i.e.\ $P_{Z_x}S_xw=0$. Moreover, for a given $w$, two admissible $y$ differ by an
element of $H_x\cap E_4(R_x)$.
\end{proof}

\begin{corollary}\label{C:pencil}
With the same notation as above,
\[
\dim\ker J(\tfrac\pi4)=\dim\ker(B_x- S_x)\quad \textrm{and}\quad\dim\ker J(\tfrac{3\pi}{4})=\dim\ker(B_x+ S_x).
\]
\end{corollary}
\begin{proof}
It follows directly from Equation \eqref{Equation-not-piover2}, since at $t=\tfrac\pi4$ and $t=\tfrac{3\pi}{4}$ we have $c=\cot t=\pm 1$, and hence
\[
c^{2}\id_W-c\,S_x-C_x=(\id_W-C_x)\mp S_x=B_x\mp S_x.\qedhere
\]
\end{proof}

\section{Polar foliations of codimension 1}\label{S:codim-1}

The main goal of this section is the classification of foliations of codimension one on the octonionic projective plane:
\begin{theorem}\label{T:codim1class}
Let $(\CaP,\fol)$ be a codimension 1 (polar) foliation. Then $\fol$ is congruent to either the foliation by concentric spheres around a point, or to the family of tubes around the canonical $\HP^2\subset\CaP$. Both are homogeneous.
\end{theorem}

Let $(\CaP,\fol)$ be a codimension 1 (polar) foliation, let $X$ denote the leaf space $X=\CaP/\fol$, and let $\sigma\colon\CaP\to X$ denote the canonical projection. By \cite[\S 5]{Tan98} we know that either $X=[0,\pi/2]$ and $\fol$ is the foliation by concentric spheres around a point, or $X=[0,\pi/4]$ and the two singular leaves $M_+=\sigma^{-1}(0)$, $M_-=\sigma^{-1}(\pi/4)$ have dimension $8$ and $11$ respectively. 
In the first case, $\mathcal{F}$ is homogeneous; actually, it is the orbit foliation of the cohomogeneity one action of a $\mathsf{Spin}(9)\subset \mathsf{F}_4$ on $\CaP$, see~\cite[p.~323]{BCO}.

From now on in this section we assume that $X=\CaP/\mathcal{F}=[0,\pi/4]$. The main goal is to prove that $M_+$ is a totally geodesic $\HP^2$, see Theorem \ref{T:main} below.
Such submanifolds are known to be congruent to each other~\cite[Theorem~1]{Wolf}. It is well known that the family of tubes around such an $\HP^2$ determines a homogeneous foliation on $\CaP$, whose other singular leaf is a minimal $\mathbb{S}^{11}$; moreover, this foliation is obtained through the cohomogeneity one action of an $\mathsf{Sp}(3)\mathsf{Sp}(1)$-subgroup of $\mathsf{F}_4$, see~\cite[Theorem~12.1.5]{BCO}. Thus, once we show that $M_+$ is a totally geodesic $\HP^2$, we will have that $\mathcal{F}$ is the orbit foliation of one of the two (up to orbit equivalence) cohomogeneity one actions on $\CaP$, thus proving Theorem \ref{T:codim1class}.

\subsection{The Lagrangian family of $M_+$-Jacobi fields}\label{S:Jac tensor}
Fix a point $p\in M_+$ and put $W:=T_pM_+$, $\nu:=\nu_pM_+$. By assumption, we have $\dim W=\dim\nu=8$. Hereafter everything will be studied at $p$, so we will omit the dependence on the point. Fixing a unit $x\in\nu$, let $\Lambda$ denote the Lagrangian family of $M_+$-Jacobi fields along $\gamma_x$.  By~\cite[Theorem~5.1]{Tan98} we know that the focal points of $M_+$ along $\gamma_x$ in $(0,\pi)$ occur at $\pi/4$, $3\pi/4$ with multiplicity $m_-=4$ (that is, when $\gamma_x$ intersects $M_-$), and at $\pi/2$ with multiplicity $7$ (when $\gamma_x$ intersects $M_+$ again). Therefore, combining Proposition \ref{prop:pencil} and Corollary \ref{C:pencil}, and using the notation and definitions introduced in \S \ref{subsec:lagrangian} applied to $M_+$, we readily obtain:

\begin{lemma}\label{L:kernels}
With the notation as above, we have:
\begin{equation}\label{Equationpiover4}
\dim\ker\bigl(B_x-S_x\bigr)=\dim\ker\bigl(B_x+S_x\bigr)=4,
\end{equation}
\begin{equation}\label{Equationpiover2}
7=\dim\bigl(H_x\cap E_4(R_x)\bigr)+\dim\ker\bigl(P_{Z_x}S_x|_{Z_x}\bigr).
\end{equation}
\end{lemma}

\begin{proposition}\label{P:four}
$\dim Z_x =\dim (W\cap E_1(R_x))\ge4$ for every unit $x\in\nu$, and if equality occurs then
$P_{Z_x}S_x|_{Z_x}=0$, i.e.\ $\langle S_xw,w'\rangle=0$ for all $w,w'\in Z_x$.
\end{proposition}
\begin{proof}
From Equation \eqref{Equationpiover2}, using that $H_x\cap E_4(R_x)=(W+ E_1(R_x))^\perp=Z_x^\perp$ in $x^\perp$, we have
\begin{align*}
7={}&15-\dim\bigl(W+ E_1(R_x)\bigr)+\dim\ker\bigl(P_{Z_x}S_x|_{Z_x}\bigr)\\
={}&15-\dim W-\dim E_1(R_x) +\dim (W\cap E_1(R_x))+\dim\ker\bigl(P_{Z_x}S_x|_{Z_x}\bigr)\\
\leq& -1+2\dim Z_x,
\end{align*}
where the last inequality uses the fact that $\ker\bigl(P_{Z_x}S_x|_{Z_x}\bigr)\subseteq Z_x$. Finally, if $\dim Z_x=4$ the inequalities above become equalities and in particular $\ker\bigl(P_{Z_x}S_x|_{Z_x}\bigr)= Z_x$, which implies $P_{Z_x}S_x|_{Z_x}=0$.
\end{proof}

\begin{proposition}\label{T:notline}
$\nu$ is not an octonionic line.
\end{proposition}
\begin{proof}
Assume for a contradiction
that $\nu$ is an octonionic line. Let $x\in\nu$ be of unit norm. Then $\nu=\ell(x)$ (by Lemma \ref{lem:cayley}) and hence $W=\ell(x)^\perp=E_1(R_x)$,
i.e.\ $\Pi_x^1=P_W$. Therefore
\[
B_x=P_W\Pi_x^1\big|_W=\id_W .
\]
By \eqref{Equationpiover4}, $\ker(\id_W-S_x)$ and $\ker(\id_W+S_x)$ both have dimension $4$; as $\dim W=8$ and $S_x$
is symmetric, the spectrum\footnote{Here we denote the multiplicity of an eigenvalue as an exponent. Thus $1$ and $-1$ are eigenvalues with multiplicity $4$.}
of $S_x$ is $\spec S_x=\{1^{4},(-1)^{4}\}$ and so $S_x^{2}=\id_W$.  By homogeneity
$S_x^{2}=|x|^{2}\id_W$ for every $x\in\nu$, and since $x\mapsto S_x$ is linear, polarisation gives
\[
S_xS_y+S_yS_x=2\langle x,y\rangle \id_W,\qquad \text{for any } x,y\in\nu\cong\RR^{8}.
\]
Thus $W$ is an 8-dimensional module over the Clifford algebra
$\mathrm{Cl}_{0,8}\cong M_{16}(\RR)$,
which contradicts the fact that such modules have dimension $16k$ for some $k\in \mathbb{N}$, see~\cite[Ch.~I, Theorem~4.3]{LM}.
\end{proof}

\begin{lemma}\label{L:4dim}
Every octonionic line $\ell\subset T_p\CaP$ intersects $\nu$ in a $0$- or $4$-dimensional subspace. If $\ell, \ell'$ are distinct octonionic lines intersecting $\nu$ nontrivially, then $\nu=(\nu\cap\ell)\oplus(\nu\cap\ell')$.
\end{lemma}
\begin{proof}
Let the octonionic line $\ell\subset T_p\CaP$ intersect $\nu$ nontrivially, and let $x\in\nu\cap\ell$ be a unit vector. Then $\ell=\ell(x)$ and, by Proposition~\ref{P:four},
\begin{align*}
\dim(\nu\cap \ell(x))={}&16-\dim(W+\ell(x)^\perp)=16-\dim(W+E_1(R_x))\\
={}&16-(\dim W+\dim E_1(R_x)-\dim(W\cap E_1(R_x)))\geq 4.
\end{align*}

By Proposition~\ref{T:notline}, there exists an octonionic line $\ell'\neq\ell$ also intersecting $\nu$ nontrivially. For any such $\ell'$, the two subspaces $\nu\cap \ell$, $\nu\cap \ell'$ intersect trivially (by Lemma \ref{lem:cayley}), have dimension $\ge 4$, and are both contained in the 8-dimensional space $\nu$. Therefore
\[
\dim(\nu\cap\ell)=\dim(\nu\cap \ell')=4,\qquad \nu=(\nu\cap \ell)\oplus(\nu\cap \ell').\qedhere
\]
\end{proof}

\subsection{Normal form}\label{S:Normal-form}
We will now derive a result that provides a convenient description of the normal space $\nu=\nu_pM_+$ as a subspace of $T_p\CaP\simeq \mathbb{O}^2$.

\begin{proposition}\label{P:normal}
Let $x\in \nu$, $|x|=1$. Up to identifying $T_p\CaP$ with $\OO^2$,
there are $4$-planes $E,D\subset\OO$ and $m_0\in\OO$ with
\[
\nu=\{(e,\,m_0e+d)\ :\ e\in E,\ d\in D\},
\]
and $x=(0,1)\in\OO^2$. 
Furthermore, considering the real vector subspace $M=\{m\in\OO: m\cdot E\subseteq D\}$ of $\OO$, we have that $\dim M=4$, one can assume $m_0\in M^\perp$, and the multiplication
\[
M\times E\longrightarrow D,\qquad (m,e)\mapsto me,
\]
is a bilinear map between $4$-dimensional Euclidean spaces with $|me|=|m|\,|e|$.
\end{proposition}
\begin{proof}
Fix a unit $x\in \nu$. Then $\nu\cap \ell(x)$
is 4-dimensional by Lemma \ref{L:4dim}. Since $\mathsf{Spin}(9)$ is transitive on unit vectors, we can identify $T_p\CaP\simeq\mathbb{O}\oplus\mathbb{O}$ by a map that sends $x$ to $(0,1)$, $\ell(x)$ to $0\oplus\mathbb{O}$ and $\ell(x)^\perp$ to $\mathbb{O}\oplus 0$. Let $\pi_i\colon\OO^2\to \OO$, $i=1,2$, denote the projections onto the first and second coordinates, and let $E:=\pi_1(\nu)$. Then $\ker\pi_1|_\nu=\nu\cap \ell(x)$ is $4$-dimensional by Lemma \ref{L:4dim} (and isomorphic to its own projection $\pi_2(\nu\cap \ell(x))=:D$), hence $E$ is 4-dimensional as well.


Take $m_0\in\OO$ such that the octonionic line $\ell':=\{(q,m_0 q): q\in \OO\}$ intersects $\nu$ nontrivially (which exists by Proposition \ref{T:notline}). By Lemma \ref{L:4dim}, 
\[\nu=(\nu\cap\ell(x))\oplus(\nu\cap\ell')=(0\oplus D )\oplus (\nu\cap\ell'),\]
and both summands are $4$-dimensional. The second summand is clearly contained in the $4$-dimensional space $\{(e, m_0 e) : e\in E\}$, hence they are equal. Therefore
\[
\nu=\{(e,m_0e+d)\ :\ e\in E,\ d\in D\}.
\]

Since the previous equation remains valid if we replace $m_0$ with $m_0+m$ for any $m\in M$, we may assume $m_0\in M^\perp$.

It remains to show $\dim M=4$. Fix $e_*\in E\setminus\{0\}$, and consider the map $M\to D$ given by $m\mapsto me_*$. It is injective because $\OO$ has no zero divisors, so it suffices to show it is surjective. Let $d\in D$. Define $z=(e_*,m_0 e_*+d)$, and let $\hat{z}\in\OO$ be the slope of the octonionic line $\ell(z)$,  that is, the octonion such that  $\ell(z)=\{(q,\hat{z}q): q\in \OO\}$. Since $\ell(z)$ intersects $\nu$ nontrivially, the same argument as above (with $\hat{z}$ in place of $m_0$) shows that
\[
\nu=\{(e,\hat{z}e+d)\ :\ e\in E,\ d\in D\}=\{(e,m_0e+d)\ :\ e\in E,\ d\in D\}.
\]
The last equation implies $\hat{z}-m_0\in M$. Since $z\in\ell(z)$, the definition of slope gives $\hat{z} e_*=m_0 e_*+d$, that is, $(\hat{z}-m_0) e_*=d$, finishing the proof that $\dim M=4$.
\end{proof}

\subsection{The spectrum of $B_x$}
\begin{proposition}\label{P:Ca}
Let $x\in\nu$, $|x|=1$, be arbitrary. Using the identification $T_p\CaP\simeq \mathbb{O}^2$ described in Proposition~\ref{P:normal} associated to $x$, we have:
\[
W=\{(-\overline{m_0}w+u,\,w)\ :\ u\in E^\perp,w\in D^{\perp}\},
\]
and the spectrum of $B_x$ is
\[
\spec B_x=\Bigl\{1^{4},\ \Bigl(\tfrac{|m_0|^2}{1+|m_0|^{2}}\Bigr)^{4}\Bigr\}.
\]
\end{proposition}
\begin{proof}
Fix a unit normal $x\in \nu$ throughout the proof. If $v, w\in \OO$ are such that $(v,w)\in W=\nu^\perp$ then for every $e\in E$ and $d\in D$ we have
\[
0=\langle(v,w), (e, m_0e+d)\rangle=\langle v,e\rangle+\langle w, m_0e+d\rangle.
\]
Letting $e=0$ we get $w\in D^\perp$, and the equation becomes
\[
0=\langle v,e\rangle+\langle w, m_0e\rangle=\langle v,e\rangle+\langle \overline{m_0}w, e\rangle=\langle v+\overline{m_0}w, e\rangle \quad \text{for all } e\in E.
\]
Hence $u=v+\overline{m_0}w\in E^\perp$, and thus 
$v=-\overline{m_0}w+u$ for $u\in E^\perp$, which by dimension reasons confirms the description of $W$. Write it as
\[
W=W_1\oplus W_2,\qquad W_1=\{(-\overline{m_0}w,w): w\in D^\perp\},\quad W_2=\{(u,0): u\in E^\perp\}.
\]
We claim that the splitting is orthogonal. First note that, for any nonzero $m\in M$, $L_m\colon \mathbb{O}\to \mathbb{O}$ is proportional to an orthogonal map (by the octonionic relation $\langle ma,mb\rangle=|m|^2\langle a, b\rangle$), and hence $L_m(E)=D$ implies $L_m(E^\perp)=D^\perp$. Let $u\in E^\perp\setminus\{0\}$ and $w\in D^\perp$ be arbitrary. Since $Mu=D^\perp$, there exists $m\in M$ such that $mu=w$. Then $\langle  w, m_0u\rangle=\langle mu, m_0 u\rangle=|u|^2\langle m,m_0\rangle=0$, since $m_0\in M^\perp$. By dimension reasons it follows that $m_0E^\perp\subseteq D^{\perp\perp}=D$. Therefore $\langle (-\overline{m_0}w,w),(u,0)\rangle=\langle  w, m_0u\rangle=0$, which proves orthogonality.

Recall that $B_x=P_W\Pi_x^1\vert_W=P_W P_{\ell(x)^\perp}\vert_W$. Under the identification $T_p\CaP\simeq \mathbb{O}^2$ associated to $x$, we have $E_1(R_x)=\ell(x)^\perp\simeq \mathbb{O}\oplus 0$ and $E_4(R_x)=\ell(x)\cap x^\perp\simeq 0\oplus\mathbb{O}$. Hence, $B_x\vert_{W_2}=\id_{W_2}$. Let us now analyse $B_x\vert_{W_1}$. Observe that for any $w\in D^\perp$ we have
\begin{equation}\label{E:W_nu_projections}
P_{\ell(x)^\perp}(-\overline{m_0}w,w)=(-\overline{m_0}w,0)=\tfrac{|m_0|^2}{1+|m_0|^{2}}(-\overline{m_0}w,w)-\tfrac{1}{1+|m_0|^{2}}(\overline{m_0}w,|m_0|^2w).
\end{equation}
Let us argue that the right-hand side is the decomposition into $W$ and $\nu$-components. Clearly, the first summand lies in $W$. To show that the second lies in $\nu$ it is enough to argue that $\overline{m_0}w\in E$. But this follows from the fact $m_0E^\perp=D$ shown above, since for every $u\in E^\perp$ we have $\langle \overline{m_0} w,u\rangle=-\langle w,m_0u\rangle=0$.
Therefore, by Equation~\eqref{E:W_nu_projections} we deduce $B_x\vert_{W_1}=P_WP_{\ell(x)^\perp}\vert_{W_1}=\tfrac{|m_0|^2}{1+|m_0|^{2}}\id_{W_1}$, which together with $B_x\vert_{W_2}=\id_{W_2}$ proves the result.\qedhere
\end{proof}

\subsection{$M_{+}$ is a totally geodesic $\HP^2$}
\begin{theorem}\label{T:main}
$M_+$ is congruent to the canonical $\HP^2\subset\CaP$.
\end{theorem}

\begin{proof}
We first prove that $M_{+}$ is totally geodesic. Fix a point $p\in M_+$ and a unit $x\in \nu=\nu_pM_+$. 
Here we will use information about $B_x$ and the fact that $\dim(\ker B_x\pm S_x)=4$ to obtain that $S_x=0$.

Define $\beta=\beta(x):=\frac{|m_0|^{2}}{1+|m_0|^{2}}$, so that $B_x$ has eigenvalues $1$ and $\beta$, both with multiplicity $4$. Using the splitting $W=Z_x\oplus Z_x^\perp$, where we recall that $Z_x=\ker C_x=\ker (\id_W-B_x)$, we can write a vector in $W$ as $(z,z')$ with $z\in Z_x$ and $z'\in Z_x^\perp$, so that $B_x(z,z')=(z,\beta z')$. Furthermore, since $\dim Z_x=\dim \ker(\id_W-B_x)=4$ (by Proposition \ref{P:Ca} has eigenvalue 1 with multiplicity exactly 4 since $\beta<1$), by the equality case of Proposition \ref{P:four} we know that $S_x(Z_x)\subseteq Z_x^\perp$, and therefore $S_x(z,z')=(X^*z',Xz+Yz')$ for some linear maps $X\colon Z_x\to Z_x^\perp$ (and $X^*$ its adjoint) and $Y\colon Z_x^\perp\to Z_x^\perp$ self-adjoint.

Let us compute $\ker(B_x\pm S_x)$:
\begin{align*}
\ker(B_x\pm S_x)={}&\{(z,z'): (z\pm X^*z',\beta z'\pm(Xz+Yz'))=(0,0)\}\\
={}&\{(z,z'): z=\mp X^*z', \beta z'-XX^*z'\pm Yz'=0\}\\
={}&\{(\mp X^*z',z'): \beta z'-XX^*z'\pm Yz'=0\}\\
\simeq{}&\ker(\beta \id_{Z_x^\perp}-XX^*\pm Y).
\end{align*}
On the other hand, since $\beta \id_{Z_x^\perp}-XX^*\pm Y$ is an operator on a 4-dimensional space, and its kernel is $4$-dimensional (because $\dim \ker(B_x\pm S_x)=4$ by Lemma~\ref{L:kernels}), it follows that
\[
\beta \id_{Z_x^\perp}-XX^*+ Y=\beta \id_{Z_x^\perp}-XX^*- Y=0,
\]
thus giving $Y=0$ and $\beta \id_{Z_x^\perp}=XX^*$.

We will get to the conclusion
\todo[Mar]{Modified the rest of the proof, there was a gap before}
by first showing that $\beta(x_*)=0$ for some unit\todo[Mig]{Added ``unit" in a couple of places} $x_*\in \nu$, and then by showing that in this case one must have $S\equiv 0$.

Suppose by contradiction that $\beta(x)\neq 0$ for all unit $x$. Then from the second equation ${1\over \sqrt{\beta}}X$ is orthogonal, from the first equation $S_x(z,z')=(X^*z',Xz)$, and then
\[
S_x^2(z,z')=(X^*Xz,XX^*z')=\beta(z,z'),\qquad\text{and hence}\qquad S_x^2=\beta \id_W.
\]
But again by polarisation (cf.\ the proof of Proposition~\ref{T:notline}), this implies that $\left\{{1\over \sqrt{\beta}}S_x: x\in \nu\right\}$ would determine an 8-dimensional Clifford system in $W\simeq \RR^8$, contradicting the fact that the lowest-dimensional nontrivial Clifford module over $\mathrm{Cl}_{0,8}$ is of dimension $16$.

By the contradiction, there exists some $x_*\in\nu$ such that $\beta(x_*)=0$. Then $m_0=0$ in the normal form associated to $x_*$, so
\[
\nu=(0\oplus D)\oplus(E\oplus0)=\bigl(\nu\cap\ell(x_*)\bigr)\oplus\bigl(\nu\cap\ell(x_*)^\perp\bigr).
\]
If $x\in\nu\cap\ell(x_*)$ is a unit vector, then $\ell(x)=\ell(x_*)$, so $\beta(x)=\beta(x_*)=0$. If $y\in\nu\cap\ell(x_*)^\perp$ is a unit vector, then $\ell(y)=\ell(x_*)^\perp$, and hence
\[
B_{y}=P_WP_{\ell(x_*)}|_W=C_{x_*}=\id_W-B_{x_*},
\]
whose spectrum is $\{0^4,1^4\}$. Therefore $\beta(y)=0$. However, since $S_x$ is symmetric and $S_x^2=\beta(x)\id_W$, it follows that $S_x=0$ for any $x\in \nu\cap \ell(x_*)$ and $S_y=0$ for any $y\in \nu\cap \ell(x_*)^\perp$. Since $\nu=(\nu\cap \ell(x_*))\oplus (\nu\cap \ell(x_*)^\perp)$ and $x\mapsto S_x$ is linear, $S_x=0$ for every $x\in\nu$.

Hence, $M_+$ is a closed $8$-dimensional totally geodesic submanifold of $\CaP$.
By Wolf's classification~\cite[Theorem 1(5)]{Wolf}, $M_+$ is either a totally geodesic $\mathbb{OP}^1$ or a totally geodesic $\mathbb{H P}^2$. The first case is impossible, as we showed in Proposition~\ref{T:notline}  that $\nu$ (and hence $M_+$) is not octonionic. This concludes the proof.
\end{proof}

\section{Polar foliations of codimension 2}\label{S:codim-2}
In this section we deal with the codimension $2$ case. The main result is:
\begin{theorem}\label{T:codim2class}
	Let $(\CaP,\fol)$ be a codimension 2 polar foliation. Then $\fol$ is either homogeneous or congruent to the inhomogeneous foliation described in Proposition~\ref{P:nonhom-fol} below.
\end{theorem}

The proof consists in first understanding the geometry of the leaf space $X=\CaP/\fol$, more precisely proving that it is a spherical triangle with very specific angles (Proposition~\ref{P:triangles}), and then analysing each case separately (Proposition~\ref{prop:Spin8}, Theorem \ref{T:pi/4} and Theorem~\ref{T:irred}).

\begin{proposition}\label{P:triangles}
Let $(\CaP,\fol)$ be a codimension 2 polar foliation. Let $X=\CaP/\fol$ denote the leaf space. Then $X$ is isometric to a spherical triangle with angles $(\pi/2, \pi/2, \pi/g)$ for $g\in\{2,4,6\}$ or $(\pi/2,\pi/3,\pi/4)$.
\end{proposition}
\begin{proof}
By \cite[Lemma 2.3]{Lyt14}, since $\CaP$ is irreducible, the foliation $\fol$ cannot be regular. By \cite[Proposition 9.2]{Lyt14} the section is a round $\RP^2$ with sectional curvature equal to 1. Since $\CaP$ is simply connected,  \cite[Theorem 1.6 and Corollary 1.7]{Lyt10} imply that $\fol$ has no exceptional leaves, and the singular leaves correspond exactly to the boundary of $X$.

Moreover, $X$ is isometric to the quotient $\mathbb{S}^2/\Gamma$, where $\Gamma$ is a finite group generated by reflections in hyperplanes (that is, great circles), see \cite[Section~2.1]{Lyt14}. Thus $\Gamma$ is a subgroup of $\mathsf{O}(3)$ generated by reflections. Since the section is $\RP^2$, the antipodal map belongs to $\Gamma$. Thus, if $\Gamma$ is reducible, it must then be of the form $\{\pm 1\}\times \mathsf{D}_{2a}$ for $a$ even, where $\mathsf{D}_{2a}$ denotes the dihedral group with $2a$ elements. If $\Gamma$ is irreducible, it must be the symmetry group of the cube or dodecahedron (because the tetrahedron is not centrally symmetric). Therefore $X$ is isometric to a spherical triangle with angles $(\pi/2, \pi/2, \pi/a)$ for some even $a$, to $(\pi/2,\pi/3,\pi/4)$, or to $(\pi/2,\pi/3,\pi/5)$.

On the other hand, recall from Section \ref{S:inf-fol} that if $L_p$ denotes a corner leaf, the corresponding infinitesimal foliation $(\nu_pL_p,\fol_p)$, when restricted to the unit sphere $\sph_p$, yields a singular foliation  $(\sph_p,\fol_p)$ of codimension 1 with singular leaves (corresponding to the directions of the singular edges). This is an isoparametric foliation of the sphere $\sph_p$, and thus, by \cite[Satz A]{Mun81},   the leaf space $\sph_p/\fol_p$ must be isometric to $[0,\pi/g]$ with $g\in \{1,2,3,4,6\}$.

The space of directions of $X$ at $\sigma(p)$ is given by the quotient of $\sph_p/\fol_p$ by a finite group acting on it by isometries (the holonomy group of the leaf $L_p$, see e.g.~\cite[Sections 3.1, 3.3]{MR19}). But the only nontrivial isometry of $[0,\pi/g]$ fixes the center of the interval. This point would correspond to a direction $x\in \sph_p$ such that for $\epsilon>0$ small, $L_{\exp_p(\epsilon x)}$ is an exceptional leaf, contradicting the non-existence of such leaves.

Thus the angle at each corner of the spherical triangle $X$ must be $\pi/g$ for some $g\in \{2,3,4,6\}$ ($g=1$ is excluded since otherwise $X$ would not have a corner at $\sigma(p)$), thus finishing the proof.
\end{proof}

Given $(\CaP,\fol)$ with leaf space $X$ isometric to a spherical triangle, label the three edges with $e_i$, and for $i=1,2,3$ define the \emph{multiplicity}, $m_i$, of $e_i$, as the difference $\dim \fol-\dim L_{p_i}$ where $p_i$ is a point projecting to $e_i$. Since $\dim \fol=14$ we can rewrite $m_i=14-\dim L_{p_i}$. We furthermore call $e_{ij}=e_i\cap e_j$, and define its multiplicity $m_{ij}$ similarly.

\subsection{The triangles $(\pi/2,\pi/2,\pi/g)$ with $g\in\{2,4,6\}$}\label{subsec:22g}

In the case of $X$ being a triangle with angles $(\pi/2,\pi/2,\pi/g)$, this is an isosceles triangle, so let us call $e_1, e_2$ the two equal edges and $e_3$ the last one. In this case, $X$ admits a submetry $X\to [0,\pi/2]$, hence $\fol$ is contained in a codimension 1 foliation with quotient $[0,\pi/2]$. By \cite[\S 5]{Tan98} this must\todo[Ric]{added citation} be the foliation by concentric spheres around a point $p_0$, thus $\{p_0\}$ is itself a leaf of $\fol$ mapping to $e_{12}$. Moreover, the cut locus $N$ of $p_0$ is a totally geodesic sphere $\sph^8({1\over 2})$ with curvature 4,  $N$ is itself foliated by leaves of $\fol$, and $N/\fol$ coincides with $e_3$. In particular $(N,\fol|_{N})$ is a foliation of codimension 1. Let $g'$ be the integer such that $N/\fol|_N\simeq [0,{\pi\over g'}]$. Rescaling $N$ by a factor of two gives a foliation in $\sph^8$, with curvature 1, and quotient $\sph^8/\fol|_N\simeq \left[0,{2\pi\over g'}\right]=\left[0,{\pi\over (g'/2)}\right]$. Equation \eqref{E:relation} applied to $(\sph^8, \fol|_{N})$ gives that
\[
2\cdot 7=(g'/2)(m_++m_-),\qquad g'/2\in \{1,2,3,4,6\},
\]
where $m_\pm$ denote the multiplicities of the singular leaves of the codimension one foliation $\fol\vert_N$ of $\mathbb{S}^8$. Hence, the only possibilities are either $g'/2=2$ (i.e.\ $g'=4$) or $g'/2=1$ (i.e.\ $g'=2$).

Let $\fol_{p_0}$ be the slice foliation at the leaf $L_{p_0}=\{p_0\}$ restricted to the unit sphere $\sph_{p_0}= T_{p_0}^1 \CaP$. Now, the exponential map
\[
(\sph^{15}=T_{p_0}^1\CaP, \fol_{p_0})\to (N,\fol|_{N})
\]
is a Riemannian submersion (and in fact, the octonionic Hopf fibration) sending leaves to leaves, and by a dimension reason it follows that the leaves of $\fol_{p_0}$ are exactly the preimages of the leaves of $\fol|_{N}$. In particular
\[
\sph^{15}/\fol_{p_0}\simeq {N}/\fol|_{N}\simeq [0,\pi/g'].
\]
On the other hand, $\sph^{15}/\fol_{p_0}$ is the space of directions of $X$ at $e_{12}$. This implies $g=g'$, and therefore $X$ is a triangle with angles either $(\pi/2, \pi/2,\pi/2)$ or $(\pi/2, \pi/2,\pi/4)$. We analyse each one of these two cases in the following two results.

\begin{proposition}\label{prop:Spin8}
There is, up to congruence, one polar foliation $(\CaP, \fol)$ with quotient $X$ a triangle with angles $(\pi/2,\pi/2,\pi/2)$. It is given by the orbits of a Lie subgroup $\Spin(8)\subset \Spin(9)\subset \mathsf{F}_4$.
\end{proposition}
\begin{proof}
By symmetry, all leaves projecting to the vertices of $X$ are points, call them $p_0, p_1, p_2$. Furthermore, the strata $N_i$ projecting to the edges opposite to $p_i$ are, by equidistance, the distance spheres of radius $\pi/2$ around $p_i$. In $\CaP$, these correspond to the (totally geodesic) cut loci $N_i\simeq \sph^8(4)$ of $p_i$.

Fixing $p_0$, the infinitesimal foliation $(\sph^{15}=T_{p_0}^1\CaP, \fol_{p_0})$ is an isoparametric foliation whose singular leaves are $T_{p_0}^1N_1$, $T_{p_0}^1N_2$, hence 7-dimensional spheres contained in perpendicular octonionic lines.

Since the $\Spin(9)$-action on $T_{p_0}\CaP\simeq \OO\oplus \OO$ is transitive on octonionic lines, in particular it can be used to move $T_{p_0}N_1$ to $\OO\oplus 0$ and $T_{p_0}N_2$ to $0\oplus \OO$. The foliation is thus fixed ($g=2$ foliations of $\sph^{15}$ are determined by their focal spheres, and $\fol$ is obtained by exponentiating $\fol_{p_0}$) and uniqueness follows. The existence of one such foliation is well known in the literature as the homogeneous foliation induced by $\Spin(8)$ (see for example \cite[Table~3]{GK16}).
\end{proof}

\begin{theorem}\label{T:pi/4}
Let $(\CaP,\fol)$ be a polar foliation with leaf space $X$ isometric to a spherical triangle of angles $(\pi/2,\pi/2,\pi/4)$, and let $m_1, m_2$ denote the multiplicities of the leaves projecting to the edges $e_1, e_2$ respectively, labeled so that $m_1\leq m_2$. Then $m_1+m_2=7$, and the foliation is uniquely determined by the values $(m_1, m_2)\in \{(1,6), (2,5), (3,4)\}$ up to congruence. In particular, up to congruence there are exactly three foliations of this type, the first two homogeneous, given by the actions of $\Spin(7) \cdot \Spin(2)$ and of $\Spin(6) \cdot \Spin(3)$, and the last one inhomogeneous.
\end{theorem}
\begin{proof}
Recall that $\fol$ is contained in the codimension 1 foliation of $\CaP$ whose leaves are $\sigma^{-1}(e_{12})=\{p_0\}$, the geodesic spheres at $p_0$ and the cut locus of $p_0$. The infinitesimal foliation $(\sph^{15}, \fol_{p_0})$ is an isoparametric foliation with $g=4$ distinct principal curvatures, since $\sph^{15}/\fol_{p_0}\simeq [0,\pi/g]$ is the space of directions of $X$ at $e_{12}$. If we let  $m_1, m_2$ denote the corresponding multiplicities, Equation \eqref{E:relation} implies $2\cdot 14=4(m_1+m_2)$, hence $m_1+m_2=7$ and thus $(m_1,m_2)\in\{(1,6), (2,5), (3,4)\}$.

We now show that, for each choice of multiplicities, the foliation $\fol$ is unique up to congruence in $\CaP$ (but see Remark~\ref{R:(3,4)} regarding the foliation $\fol_{p_0}$). Suppose therefore that $\fol$, $\fol'$ are two polar foliations on $\CaP$, with respective quotients $X$, $X'$, both isometric to a spherical triangle of angles $(\pi/2,\pi/2,\pi/4)$, with the same multiplicities $(m_1,m_2)\in\{(1,6), (2,5), (3,4)\}$. By an element of $\mathsf{F}_4$, we can assume that $p_0$ and $p_0'$ coincide. Let $\Spin(9)$ denote the isotropy at $p_0=p_0'$. 
The infinitesimal foliations $\fol_{p_0}$, $\fol'_{p_0}$ on $\sph^{15}=T^1_{p_0}\CaP$ factor through the foliations $(N, \fol|_N)$ and $(N, \fol'|_N)$, via the $\Spin(9)$-equivariant Hopf fibration $\sph^{15}\to N\simeq \sph^{8}(4)$, where $\Spin(9)$ acts on $N$ via $\mathsf{SO}(9)$. Now, the foliations $\fol|_N$, $\fol'|_N$ are isoparametric foliations with $2$ distinct principal curvatures and the same multiplicities $m_1, m_2$. In particular, they are congruent, and up to an element of $\mathsf{SO}(9)$ we can assume that $\fol|_N=\fol'|_N$. Therefore, up to an element of $\Spin(9)$, we can assume $\fol_{p_0}=\fol'_{p_0}$. Equality $\fol=\fol'$ follows from the fact that every leaf of $\fol$ resp.\ $\fol'$ is the exponential image of a leaf of $\fol_{p_0}$ resp.~$\fol'_{p_0}$.

This concludes the proof in the cases of multiplicities $(1,6)$ and $(2,5)$, since there exists a homogeneous polar foliation with the given multiplicities, corresponding to the action of $\Spin(7)\cdot \Spin(2)\subset \Spin(9)\subset \mathsf{F}_4$ and of $\Spin(6)\cdot \Spin(3)\subset \Spin(9)\subset \mathsf{F}_4$, respectively, see~\cite[Table~3]{GK16}.

We are left to show existence of a polar foliation on $\CaP$ with $X$ a spherical triangle of angles $(\pi/2, \pi/2, \pi/4)$ and multiplicities $(m_1, m_2)=(3,4)$. Its infinitesimal foliation at the leaf corresponding to the $\pi/4$-vertex
\todo[Ric]{It seems we don't use this fact about composed Clifford foliations, right? Should we delete this sentence, or at least have it in a separate remark after Proposition \ref{P:nonhom-fol}?}
is a special example of a composed Clifford foliation from \cite{Rad14}.
We prove this in Proposition \ref{P:nonhom-fol} below.
\end{proof} 

\begin{proposition}\label{P:nonhom-fol}
There exists an inhomogeneous codimension 2 polar foliation $(\CaP, \fol_c)$ whose quotient space is a spherical triangle with angles $(\pi/2,\pi/2,\pi/4)$ and whose multiplicities are $(m_1,m_2)=(3,4)$.
\end{proposition}
\begin{proof}
Fix a point $p_0\in \CaP$ and its cut locus $N\simeq \sph^8(4)$. On (the rescaled, constant curvature $1$) $N$, consider an isoparametric foliation $\fol_N$ with $g=2$ principal curvatures and multiplicities $(3,4)$, which is unique up to congruence in $\mathbb{S}^8$ and given by the orbits of the reducible $\mathsf{SO}(4)\times \mathsf{SO}(5)$ representation on $\RR^9=\RR^4\oplus \RR^5$. The exponential map $v\mapsto \exp_{p_0}(\tfrac\pi2 v)$
defines a Hopf fibration $\sph^{15}=T^1_{p_0}\CaP\to N$, which we can use to pull back the foliation $\fol_N$ to an isoparametric foliation $(\sph^{15},\fol_{p_0})$ with $g=4$, and quotient map
\[
\tau\colon \sph^{15}\to [0,\tfrac\pi4].
\]
($\fol_{p_0}$ is the \emph{composed foliation} $\fol_N\circ \fol_{C_{8,1}}$ using the terminology in \cite{GR16}, where it was proved to be inhomogeneous.) Finally, define $(\CaP,\fol)$ as the foliation obtained by exponentiating the leaves of $\fol_{p_0}$, i.e.\ by defining the leaf through any point $q=\exp_{p_0}(tv)$ as $\exp_{p_0}(t L_v)$, where $t\in [0,\tfrac\pi 2]$,
$v\in T^1_{p_0}\CaP$, and $L_v$ is the leaf of $\fol_{p_0}$ through $v$.

We can alternatively define the foliation above via the map:
\[
\varphi:\CaP\to \sph^2,\qquad \varphi(\exp_{p_0}tv)=\cos(t)e_3+\sin(t)(\cos(\tau(v))e_1+\sin(\tau(v))e_2),\qquad
\]
where $v\in \sph^{15}$, $t\in [0,\tfrac\pi2]$, and $\{e_1,e_2,e_3\}$ denotes here the canonical basis of $\mathbb{R}^3$\todo[Mig]{Added this, because the $e_i$ where the edges of triangles}. The image $X'=\varphi(\CaP)$ is a spherical triangle in $\sph^2$ with the north pole as vertex of angle $\tfrac\pi4$, and two vertices on the equator. The preimages of the boundary components consist of $\exp(t\cdot M_-)$, $\exp(t\cdot M_+)$ (where $M_\pm\subset \sph^{15}$ are the singular leaves of $\fol_{p_0}$ and $t\in [0,\tfrac\pi2]$) and $N$. Notice furthermore that, since $\tau$ is constant on Hopf fibres (because $\fol_{p_0}$ is pulled back from $\fol_N$) the map $\varphi$ is still well-defined at $t=\tfrac\pi2$.

We need to prove that this is a polar foliation with the desired properties.

{\it Step 0: \todo[Mig]{Changed bold to italics: we have many subsections in bold, and might be confusing to use here bold as well, inside a proof}for every $x,v\in\sph^{15}$ with $v$ horizontal vector for $L_x\subset \sph^{15}$, the plane spanned by $x, v$ has curvature 1.} In fact, recall that $E_4(R_x)=\ell(x)\cap x^\perp$ is the tangent space at $x$ of the Hopf fibre $\ell(x)\cap\sph^{15}$, hence vertical with respect to the map $\sph^{15}\to N$. Therefore, the horizontal vector $v$ lives in $\ell(x)^\perp$, and then $\operatorname{sec}(x,v)=1$ by Lemma \ref{lem:cayley}~(1).

{\it Step 1: $\fol$ is a singular foliation}, i.e.\ it is spanned by smooth vector fields. This follows from the same proof as in \cite[Proposition~3.1]{Rad14}. In short: the foliation in the complement of $N$ is smooth because it is foliated diffeomorphic to the infinitesimal foliation on $T_{p_0}\CaP$, and in a tubular neighbourhood around $N$ it is smooth because it locally looks like a product of $(N,\fol_N)$ with the foliation $(\RR^8,\fol'')$ by concentric spheres around the origin.

{\it Step 2: The restriction of $\varphi$ to the interior of $X'$ is a polar Riemannian submersion}. Consider a point $p\in \CaP$, projecting to a point in the interior of $X'$. We can write $p=\exp_{p_0} t_0 x_0$, with $t_0\in (0,\tfrac\pi2)$ and $x_0\in \sph^{15}$ a regular point for $\fol_{p_0}$. Let $\alpha\colon (-\epsilon, \epsilon)\to \sph^{15}$ denote a unit horizontal geodesic for $\fol_{p_0}$ with $\alpha(0)=x_0$ (hence $\tau(\alpha(s))=\tau(x_0)+s$) and define
\todo[Mig]{$p_0$ instead of $p$ in this exp and the one below, right?}
\todo[Mar]{Yes}
\[
F(s,t)=\exp_{p_0}(t\alpha(s))
\]
so that
\begin{align*}
\varphi(F(s,t))&=\varphi(\exp_{p_0}(t\alpha(s)))\\
&=\cos(t)e_3+\sin(t)(\cos(\tau(x_0)+s)e_1+\sin(\tau(x_0)+s)e_2).
\end{align*}
Then at $s=0$, $t=t_0$ we have
\begin{itemize}
\item $\varphi_* {\partial F \over \partial t}={\partial \varphi(F)\over \partial t}=-\sin(t_0)e_3+\cos(t_0)(\cos(\tau(x_0))e_1+\sin(\tau(x_0))e_2)$, and thus
\[
\left\|\varphi_* \tfrac{\partial F}{\partial t}\right\|^2=1.
\]
\item $\varphi_* \tfrac{\partial F}{\partial s}=\tfrac{\partial \varphi(F)}{\partial s}=\sin(t_0)(-\sin(\tau(x_0))e_1+\cos(\tau(x_0))e_2)$, thus
\[
\left\|\varphi_* \tfrac{\partial F}{\partial s}\right\|^2=\sin^2(t_0)\qquad\text{and}\qquad\left\langle \varphi_* \tfrac{\partial F}{\partial t}, \varphi_* \tfrac{\partial F}{\partial s}\right\rangle=0.
\]
\end{itemize}
At the same time, by the Gauss lemma:
\[
\|\tfrac{\partial F}{\partial t}\|^2\equiv1=\left\|\varphi_* \tfrac{\partial F}{\partial t}\right\|^2\qquad g(\tfrac{\partial F}{\partial t},\tfrac{\partial F}{\partial s})\equiv0=\left\langle \varphi_* \tfrac{\partial F}{\partial t}, \varphi_* \tfrac{\partial F}{\partial s}\right\rangle,
\]
while
\[
\|\tfrac{\partial F}{\partial s}(0,t)\|^2=\sin^2(t)=\left\|\varphi_* \tfrac{\partial F}{\partial s}\right\|^2
\]
since $J(t):=\tfrac{\partial F}{\partial s}(0,t)$\todo[Mig]{Changed from $\tfrac{\partial F}{\partial s}(\exp_{p_0}tx_0$)} is a Jacobi field with $J(0)=0$ and $\operatorname{sec}(J'(0),x_0)=\operatorname{sec}(\alpha'(0),x_0)=1$ by Step 0.

These computations imply first of all that $\varphi$ is a submersion onto the interior of $X'$, and furthermore that the restriction of $\varphi$ to the image of $F$ is an isometry. Since it is easy to see that $\|\varphi_*\|\leq 1$ in general it follows in particular that the spaces tangent to $F$ are horizontal (otherwise necessarily $\|\varphi_*\|<1$) and $\varphi$ is a Riemannian submersion.

Finally, $F$ is just the exponential map applied to $\operatorname{span}(x_0,\alpha'(0))\subset T_{p_0}\CaP$ which, by Step 0, has $\operatorname{sec}=1$. Therefore, the image $\Sigma$ of $F$ lies in a totally geodesic $\RP^2$.

{\it Step 3: $\fol$ is a (polar, inhomogeneous) singular Riemannian foliation}. We need to prove that leaves are equidistant. For leaves projecting to the interior of $X'$ this follows from Step 2, since fibers of Riemannian submersions have this property. Since the remaining leaves are limits of regular ones, equidistance of general leaves follows and thus $\fol$ is a singular Riemannian foliation. To prove that it is polar, again we use polarity of $\fol$ in the preimage of the interior of $X'$: for Riemannian submersions, this is equivalent to the vanishing of the O'Neill tensor in the regular part (e.g. by noticing that for $x,v\in T_p\Sigma$, $\sec(x,v)=1=\sec(T_p\sph^2)$ and using the O'Neill formulas). From this, Theorem 1.6 of \cite{Ale06} implies that the foliation $\fol$ itself is polar.

Inhomogeneity follows either by the fact that $\fol_{p_0}$ is inhomogeneous  \cite[Table 2]{GR16}, or by computing the multiplicities of the leaves, and checking in the classification of polar actions on $\CaP$ (cf. \cite{PT, GK16}) that these do not match with any of the possible group actions.
\end{proof}

\begin{remark}\label{R:(3,4)}
	 We note that the proof of Theorem~\ref{T:pi/4} did not make use of the classification of isoparametric foliations on spheres (except for the classical case of $g=2$). However, at the beginning of the proof we showed that $\fol_{p_0}$ is an isoparametric foliation with $g=4$ and multiplicities $(m_1, m_2)\in \{(1,6), (2,5), (3,4)\}$ on $\mathbb{S}^{15}$. By the classification of isoparametric foliations on spheres with $g=4$, for the first two multiplicity pairs there is only one (up to congruence in $\mathsf{SO}(16)$) foliation on $\mathbb{S}^{15}$, which is homogeneous, whereas for $(m_1,m_2)=(3,4)$ there are exactly two such foliations, one of them homogeneous~\cite{Chi}, \cite[\S 4.3]{FKM}. The first part of the proof of Theorem~\ref{T:pi/4} implies in particular that the inhomogeneous isoparametric foliation with $g=4$ and $(m_1,m_2)=(3,4)$ cannot be subfoliated by the fibers of the Hopf map $\mathbb{S}^{15}\to\mathbb{S}^8$.
\end{remark}

\subsection{The triangle $(\pi/2,\pi/3,\pi/4)$: a description of the section}\label{subsec:234_description}

\begin{figure}[hbt]
\centering
\begin{subfigure}{.5\textwidth}
  \centering
  \includegraphics[width=.8\linewidth]{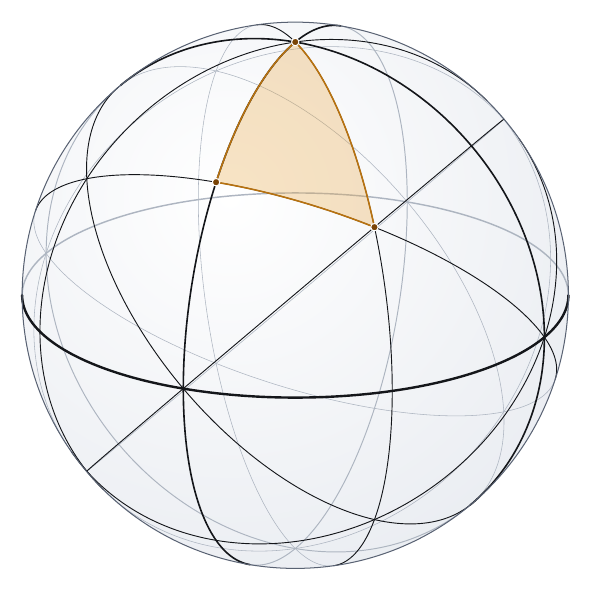}
  \caption{A fundamental domain}
  \label{fig:sub1}
\end{subfigure}%
\begin{subfigure}{.5\textwidth}
  \centering
  \includegraphics[width=.8\linewidth]{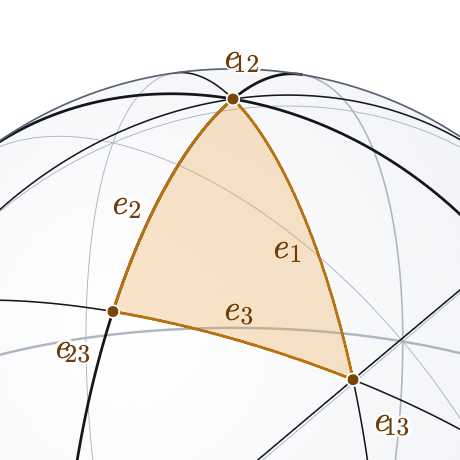}
  \caption{The labels}
  \label{fig:sub2}
\end{subfigure}
\caption{The tessellation of $\Sigma$}
\label{F:tess}
\end{figure}

Let $(\CaP,\fol)$ be now a polar foliation with leaf space $X$ isometric to a spherical triangle of angles $(\pi/2,\pi/3,\pi/4)$. Then the section $\Sigma\simeq \RP^2$ is acted on by a reflection group which, lifted to $\tilde{\Sigma}\simeq \sph^2$, is the full group of isometries of the cube. This is the group generated by the reflections in $\mathbb{R}^3$ across the planes $x_i=0$ and $x_i=\pm x_j$, for every $i\in\{1,2,3\}$ and $j\neq i$. The quotient $X$ is then isometric to a fundamental domain of this action, as shown in Figure \ref{F:tess}. We will label the edges of $X$ by letting $e_i$, $i\in\{1,2,3\}$, denote the edge opposite to the vertex with angle $\pi/(i+1)$. The vertices will be labeled $e_{ij}=e_i\cap e_j$ as above. The $e_2$ edges lie on the three great circles $x_i=0$ in $\mathbb{S}^2$, whereas the $e_1$ and $e_3$ edges lie on the circles $x_i=\pm x_j$, $j\neq i$, alternating at the $\pi/3$-vertices, see Figure~\ref{F:tess}. The leaf corresponding to the corner $e_{ij}$ (resp.\ the edge $e_i$) will be denoted $L_{ij}$ (resp.~$L_i$).

\subsection{The triangle $(\pi/2,\pi/3,\pi/4)$: computing multiplicities}

Our goal in this section is to compute the multiplicities of edges and corners.

\begin{lemma}\label{L:formula1}
Suppose $e_i$, $e_j$ have multiplicity $m_i$, $m_j$ respectively, and they intersect at an angle ${\pi/g}$. Then
\[
m_{ij}={g\over 2}(m_i+m_j).
\]
\end{lemma}
\begin{proof}
The leaf $L_{ij}$ has dimension $14-m_{ij}$, hence codimension $2+m_{ij}$, and thus $\dim(\sph_{p_{ij}})=1+m_{ij}$, for any $p_{ij}\in L_{ij}$. The foliation $\fol_{p_{ij}}$ is of codimension 1 (hence regular leaves have dimension $m_{ij}$) and the singular leaves correspond to the fibers of the projections $L_{i}\to L_{ij}$ and $L_{j}\to L_{ij}$. In particular, they have dimension $m_{ij}-m_{i}$ and $m_{ij}-m_j$. Therefore, the multiplicities of the singular leaves in $\sph_{p_{ij}}$ are $m_i$ and $m_j$. For codimension 1 foliations, Equation \eqref{E:relation} reads
\[
g(m_i+m_j)=2(\dim(\sph_{p_{ij}})-1)=2 m_{ij}
\]
and from this the formula follows.
\end{proof}

\begin{lemma}\label{L:formula2}
Let $\Sigma\simeq \RP^2$ be a section of $(\CaP,\fol)$, and suppose that a geodesic $\gamma\colon[0,\pi]\to \Sigma$ with $\gamma(0)$ regular, once projected to $X$, intersects $e_i$ for $n_i$ times, $i=1,2,3$. Then:
\[
\sum_{i=1}^3n_im_i=21
\]
\end{lemma}
\begin{proof}
Along the geodesic $\gamma$, let $c_0$ denote the number of conjugate points in $\CaP$\todo[Ric]{Added ``in $\CaP$'' to avoid confusion, since $\gamma$ is also a geodesic in $\RP^2$.}  of $\gamma(0)$ in $(0,\pi]$ ($c_0=22$ by e.g.\ Equations~\eqref{E:OperatorsGivingJ} and~\eqref{E:Jacfield} with $w=0$), and let $c_1$ denote the number of focal points of $L_{\gamma(0)}$ in $(0,\pi]$. Then, the number of conjugate (resp. focal) points in $(0,k\pi]$ becomes $kc_0$ and $kc_1$ respectively. However, by \cite[Theorem~1.1]{Lyt09}, the difference $|kc_1-kc_0|$ is bounded independently of $k$, which means $c_1=c_0$.

By \cite[\S 3.2]{Lyt09} we have $c_1=1+\sum_{i=1}^3n_im_i$, where $1$ is the index of $\gamma$ projected to $X$, and the second summand counts the intersections of $\gamma$ with the singular leaves. The result then follows.
\end{proof}

In the case of $X$ being the triangle of angles $\pi/2, \pi/3, \pi/4$ we obtain $m_1=m_3\in \{1,2,4,8\}$, because when $g=3$ the multiplicities are equal and only allowed to be the values listed here. By looking at the corresponding tessellation of the 2-sphere in Figure \ref{F:tess}, we easily compute $n_1+n_3=6$ and $n_2=3$. Therefore the formula in Lemma \ref{L:formula2} gives
\[
6m_1+3m_2=21,\quad\text{and hence}\quad 2m_1+m_2=7.
\]
We thus get that $m_1=1$ resp.\ 2, and $m_2=5$ resp.\ 3. Using Lemma \ref{L:formula1} we then have Table \ref{table:multiplicities}.\todo[Ric]{Referred to Table 2, because the Table may float around and end not exactly below this sentence.}
\begin{table}[h]
\[
\begin{array}{c||c|c|c|c|c|c}\textrm{Case}&m_1 & m_2 & m_3 & m_{12} & m_{13} & m_{23} \\\hline\hline (I)&1 & 5 & 1 & 12 & 3 & 6 \\\hline (I\!I)&2 & 3 & 2 & 10 & 6 & 5\\\end{array}
\]
\caption{Multiplicities of singular leaves}\label{table:multiplicities}
\end{table}

These are precisely the multiplicities of the remaining homogeneous examples not considered in Proposition~\ref{prop:Spin8} and Theorem~\ref{T:pi/4}, see \cite[Table~3]{GK16}. In Figure \ref{F:dim-cod} we show the dimension of leaves, depending on  the two cases:
\begin{figure}[hbt]\centering
\begin{subfigure}{.5\textwidth}
  \centering
  \includegraphics[width=.8\linewidth]{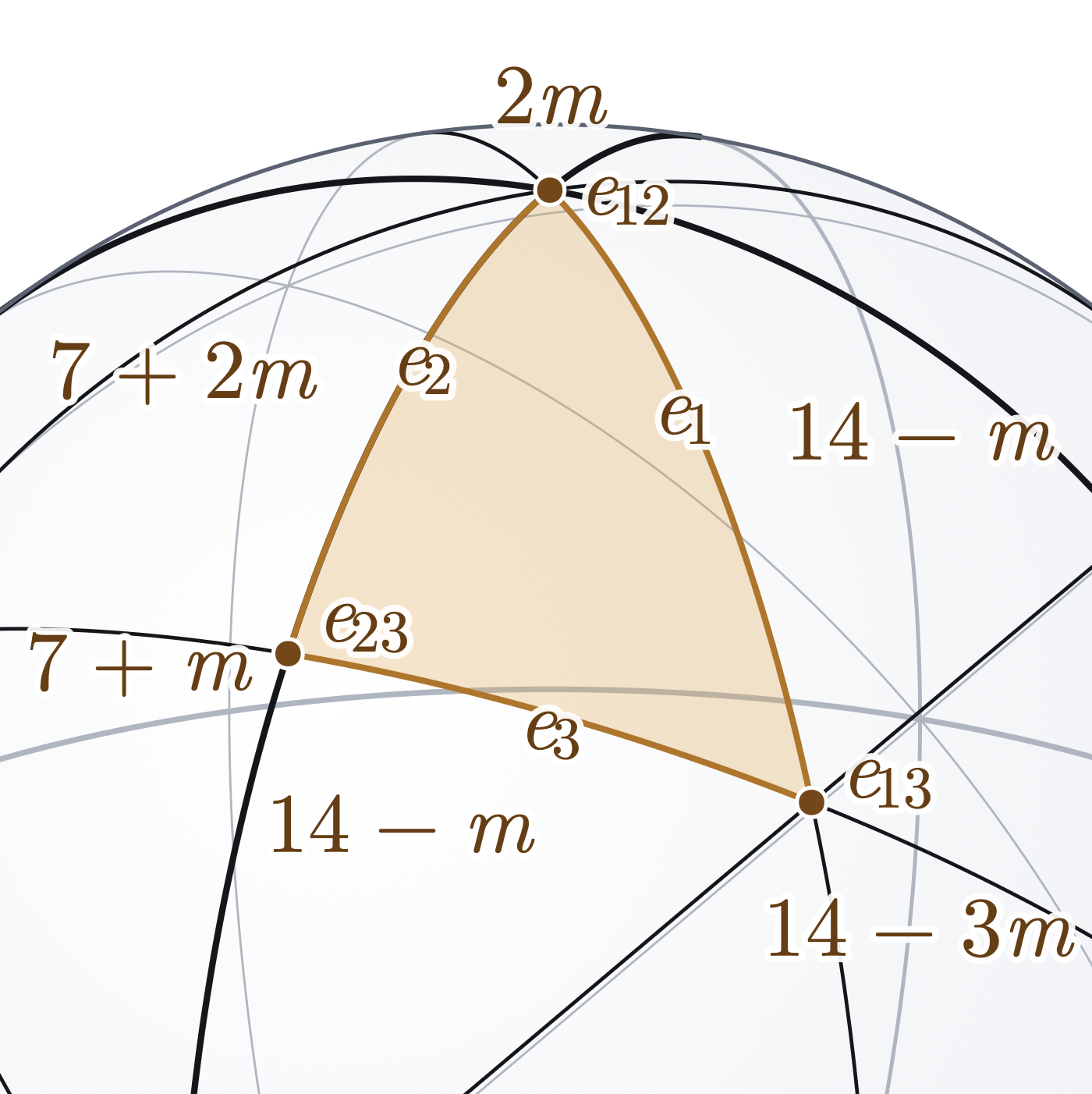}
  \caption{Dimensions}
\end{subfigure}%
\begin{subfigure}{.5\textwidth}
  \centering
  \includegraphics[width=.8\linewidth]{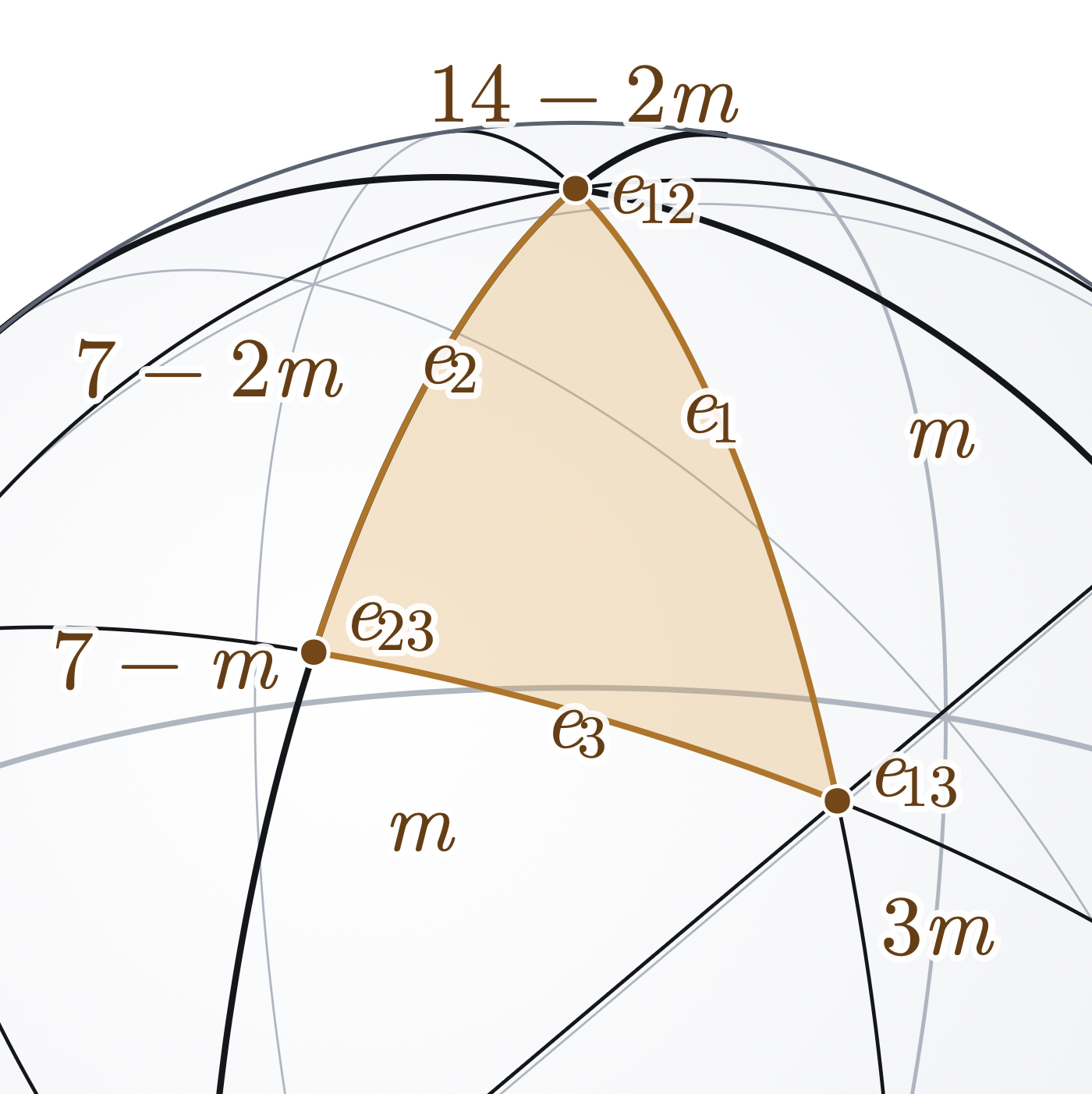}
  \caption{Multiplicities}
\end{subfigure}
\caption{Dimensions and multiplicities of singular leaves, depending on $m=m_{1}\in\{1,2\}$}
\label{F:dim-cod}
\end{figure}

\subsection{Horizontal geodesics and $L_{12}$-Jacobi fields, from $L_{12}$}

Recall that $L_{12}$ denotes the leaf corresponding to the corner $e_{12}$ with angle $\pi/4$ of the leaf space $X=\CaP/\fol$. It has dimension $2m$. Fix a point $p$ of $L_{12}$. To use the results of \S\ref{subsec:lagrangian}, we set $W:=T_pL_{12}$ and $\nu:=\nu_pL_{12}$. For every unit $x\in \nu$, let $\tau(x)\in [0,\pi/4]$
denote the angle between the projection of $\gamma_x$ in $X$, and the edge
$e_2$, where $\gamma_x(t)=\exp_p(tx)$.

In particular:
\begin{itemize}
\item For $\tau(x)=0$, the geodesic $\gamma_x$ projects to a curve entirely contained in $e_2$.
\todo[Mar]{exchanged the first two bullet points}
\item For $\tau(x)=\pi/4$, the geodesic $\gamma_x$ projects to a curve entirely contained in $e_1\cup e_3$.
\item For $\tau(x)\in(0,\pi/4)$, the geodesic $\gamma_x|_{t\in [0,\pi]}$ projects to an orbifold geodesic in $X$ passing through the interior of $X$, intersecting $e_2$ at $t=\pi/2$ (independently of $x$) with multiplicity $7-2m$, and intersecting $e_1\cup e_3$ at the four values  of $t$ given by \todo[Ric]{changed the wording slightly, ``4 times $t$'' sounded confusing to me.}
\[
\cot(t)=\pm \cos(\tau),\qquad \cot(t)=\pm \sin(\tau),
\]
all with multiplicity $m$.
\end{itemize}
From Equation \eqref{Equation-not-piover2}, we obtain that for $c=\pm \sin(\tau)$ and $c=\pm \cos(\tau)$,
\[
\dim\ker\big(c^2\id_W-cS_{x}-C_{x}\big)=m.
\]

\begin{corollary}\label{cor:traces}
For every unit vector $x\in\nu$,
\begin{equation}\label{eq:traces}
\tr S_x=0,\qquad \tr C_x+\tfrac12\|S_x\|^2=m,\qquad
\det C_x=\Big(\tfrac12\sin(2\tau)\Big)^{2m}.
\end{equation}
\end{corollary}
\begin{proof}
Start with $x$ a unit vector with $\tau(x)\in (0,\tfrac\pi4)$: this corresponds to a vector whose geodesic $\exp_ptx$ passes through regular leaves. Consider the polynomial
\[
p(c)=\det(c^2\id_W-cS_{x}-C_{x}).
\]
Since the matrix $c^2\id_W-cS_{x}-C_{x}$ has order $\dim W=2m$, $p(c)$ has degree $4m$. Since $c^2$ multiplies $\id_W$, it follows from the formula of the determinant that the coefficient of $c^{4m}$ is $\det \id_W=1$. At the same time, $p(c)$ vanishes on $\pm \sin(\tau)$, $\pm\cos(\tau)$ with multiplicity $m$, and since $\deg(p)=4m$ we have
\[
p(c)=(c^2-\sin^2(\tau))^m(c^2-\cos^2(\tau))^m=\left(c^4-c^2+{1\over 4}\sin^2(2\tau)\right)^m.
\]
However, by defining the new operator
\[
A\colon W^2\to W^2,\qquad A(w_1, w_2)=(w_2,C_xw_1+S_xw_2),
\]
it is easy to see that for any $\lambda\in \RR$, $(w_1,w_2)\in W^2$ is an eigenvector of $A$ with eigenvalue $\lambda$ if and only if $w_1\in\ker (\lambda^2\id_W-\lambda S_{x}-C_{x})$. In particular, the eigenvalues of $A$ are exactly $\pm \cos(\tau),\pm\sin(\tau)$ each with multiplicity $m$, and $p(c)$ is the characteristic polynomial of $A$. With respect to a basis of $W^2$ we can write
\[
A=\left(\begin{array}{c|c}0 & \id_W  \\\hline  C_x & S_x\end{array}\right).
\]
In particular, writing $p(c)=\sum a_ic^i$, we get:
\begin{enumerate}
\item $\operatorname{tr}(S_x)=\operatorname{tr}(A)=-a_{4m-1}=0$.
\item $\det(C_x)=\det(A)=a_0={1\over 4^m}\sin^{2m}(2\tau)$.
\item Since\todo[Mig]{Interchanged C and S in the 2,1 position}
\[
A^2=\left(\begin{array}{c|c} C_x & S_x\\\hline S_xC_x & C_x+S_x^2\end{array}\right),
\]
taking the trace gives $2\operatorname{tr} C_x+\operatorname{tr} (S_x^2)=\operatorname{tr}(A^2)$, which is equal to the sum of the squares of the roots of $p(c)$, hence
\[
2\operatorname{tr} C_x+\|S_x\|^2=2m\left(\cos^2(\tau)+\sin^2(\tau)\right)=2m.
\]
\end{enumerate}
These provide the desired results whenever $\tau(x)\in (0,\tfrac\pi4)$. By continuity of $x\mapsto S_x$ and $x\mapsto C_x$, the result follows for all $x$.
\end{proof}
\subsection{$L_{12}$ is totally geodesic}

Take $p\in L_{12}$. Given a unit vector $x\in \nu$, let $\tau=\tau(x)$ be as before.

\begin{lemma}\label{L:norm-of-S}
For any unit $x\in\nu$ with $\tau(x)\in (0,\pi/4)$,
\[
\|S_x\|^2\leq 2m(1-\sin(2\tau(x))).
\]
\end{lemma}
\begin{proof}
Since $C_x$ is positive semidefinite, we can apply the AM-GM inequality to its eigenvalues, which combined with Corollary~\ref{cor:traces} yields
\[
{1\over 2}\sin(2\tau(x))=\sqrt[2m]{\det C_x}\leq {\operatorname{tr}(C_x)\over 2m}={m-\tfrac12\|S_x\|^2\over 2m}.
\]
Reordering the terms gives the inequality.
\end{proof}

\begin{corollary}
$L_{12}$ is totally geodesic.
\end{corollary}
\begin{proof}
Fix $p\in L_{12}$ and take the infinitesimal foliation $(\sph_p,\fol_p)$. As we said, it is a codimension 1 foliation, with leaf space isometric to the interval $[0,\pi/4]$, and projection $\tau\colon \sph_p\to [0,\pi/4]$. Lemma \ref{L:norm-of-S} can be restated by saying that for $t\in (0,\pi/4)$, any $x\in\tau^{-1}(t)$ satisfies $\|S_x\|^2\leq 2m(1-\sin(2t))$. Since the map $x\mapsto \|S_x\|^2$ is continuous, it follows that for any $x\in M_-=\tau^{-1}(\pi/4)$ one has $\|S_x\|^2=0$, hence $S_x=0$. However, since  $\fol_p$ is irreducible  (i.e.\ $\operatorname{diam}(\sph_p/\fol_p)<\tfrac\pi2$) it follows that $M_-$ is \emph{full} in $\nu_pL_{12}$ (i.e.\ $\operatorname{span}(M_-)=\nu_pL_{12}$). Since $x\mapsto S_x$ is linear, it follows that $S_x=0$ for any $x\in \nu_pL_{12}$.
\end{proof}

\subsection{A characterization of $L_{12}$}
It is time to identify $L_{12}$ among the possible totally geodesic submanifolds of $\CaP$. Recall that by~\cite[Theorem~1]{Wolf},
every such manifold is congruent to the standard $\RP^2,\CP^2, \HP^2$ or $\sph^n(4)$, $1\leq n\leq8$, where $\mathbb{S}^8(4)=\mathbb{O}\mathbb{P}^1$.
\begin{proposition}We have:
\begin{enumerate}[\rm(a)]
\item If $L_{12}$ is 2-dimensional (case (I) of Table~\ref{table:multiplicities}, with $m=1$), then it is a totally geodesic $\RP^2$ of curvature 1.
\item If $L_{12}$ is 4-dimensional (case (II) of Table~\ref{table:multiplicities}, with $m=2$), then it is a totally geodesic $\CP^2$.
\end{enumerate}
\end{proposition}
\begin{proof}
(a) Assume we are in the case $m=1$. Take an orthonormal basis $\{v_1,v_2,x_1,\ldots, x_{14}\}$ of $T_p\CaP$ with $v_i\in W$ and $x_j\in \nu$. By Remark \ref{R:C-and-curv}, and using the formulas in \eqref{eq:traces}, we have
\[
1=\operatorname{tr}(C_{x_i})=\tfrac13(\sec(v_1,x_i)+\sec(v_2,x_i)-2),
\]
and hence $\sec(v_1,x_i)+\sec(v_2,x_i)=5$. Summing over all $x_i$, we get
\[
70=\sum_{i=1}^{14}\sum_{j=1}^2\sec(v_j,x_i)=\operatorname{Ric}(v_1)+\operatorname{Ric}(v_2)-2\sec(v_1,v_2)=72-2\sec(v_1,v_2),
\]
where we used the fact that $\CaP$ has constant Ricci curvature $\operatorname{Ric}(v)=36\|v\|^2$. From this, one gets the ambient curvature $\sec(v_1,v_2)=1$. However, since $L_{12}$ is totally geodesic, it follows by the Gauss equation that this is the intrinsic curvature of $L_{12}$ as well.

(b) Assume now that we are in the case $m=2$. As before, take an orthonormal basis $\{v_1,\ldots,v_4,x_1,\ldots, x_{12}\}$.  Again by Remark \ref{R:C-and-curv} and the formulas in \eqref{eq:traces}:
\[
2=\operatorname{tr}(C_{x_i})=\tfrac13\Biggl(\Big(\sum_{j=1}^4\sec(v_j,x_i)\Big) -4\Biggr),\quad \text{and then}\quad \sum_{j=1}^4\sec(v_j,x_i)=10.
\]
Summing over all $x_i$ and again using the formulas for the Ricci curvature of $\CaP$, we get:
\[
120=\sum_{i=1}^{12}\sum_{j=1}^4\sec(v_j,x_i)=\sum_{j=1}^4\operatorname{Ric}(v_j)-\sum_{i\neq j=1}^4\sec(v_i,v_j)=144-\operatorname{scal}(L_{12}),
\]
and hence  $\operatorname{scal}(L_{12})=24$.
This uniquely identifies $\CP^2$ among the possible 4-dimensional totally geodesic subspaces, since $\operatorname{scal}(\sph^4(4))=48$.
\end{proof}

We can finally conclude the analysis initiated in~\S \ref{subsec:234_description}.
\begin{theorem}\label{T:irred}
Let $(\CaP,\fol)$ denote a codimension 2 polar foliation whose quotient is the spherical triangle with angles $(\pi/2,\pi/3,\pi/4)$. Then $\fol$ is homogeneous.
\end{theorem}
\begin{remark}\label{rem:234}
	Combining Theorem~\ref{T:irred} with \cite[Table~3]{GK16} we deduce that $\fol$ is the orbit foliation of the action of subgroups $\mathsf{SU}(3)\cdot \mathsf{SU}(3)$ or $\mathsf{SO}(3)\cdot \mathsf{G}_2$ on $\CaP$.
\end{remark}
\begin{proof}
Let $m\in\{1,2\}$ be the multiplicity of the leaves projecting to the edge $e_1$ opposite to the angle $\pi/2$. In each case, there is a unique homogeneous foliation $(\CaP,\fol')$, up to congruence, with the same quotient and multiplicities, see~\cite[Table~3]{GK16}. We will denote with a ``$\,'\,$'' the object related to the homogeneous model.

Since what we obtained in this section does not depend on the specific foliation but only on the shape of the section, it follows that the leaf $L_{12}'$ is totally geodesic and isometric to $\RP^2$, resp.\ $\CP^2$, in the case of $m=1$, resp.\ $m=2$. Again by Wolf's classification~\cite[Theorem~1]{Wolf}, it follows that there is an ambient isometry moving $L_{12}$ to $L_{12}'$. Thus we can assume $L_{12}=L_{12}'$. In particular, for any $p\in L_{12}$, we have $W_p=W_p'$, $\nu_p=\nu_p'$ and, for any $x\in \nu_p$, $C_x=C_x'$ because $C_x$ only depends on $W_p$.

Look now at the infinitesimal foliations $(\sph_p,\fol_p)$, resp.\ $(\sph_{p'},\fol_p')$. These are uniquely determined by functions $\tau\colon \sph_{p'}\to [0,\pi/4]$, resp. $\tau'\colon\sph_{p'}\to [0,\pi/4]$. However, by Corollary \ref{cor:traces} we have
\[
\sin(2\tau(x))=2\sqrt[2m]{\det C_x}=2\sqrt[2m]{\det C'_x}=\sin(2\tau'(x))
\]
and since $\tau(x), \tau'(x)\in [0,\pi/4]$, we obtain that $\tau(x)=\tau'(x)$ for every $p\in L_{12}$ and every $x\in \nu_p$.

For a given $\theta\in [0,\pi/4]$ and a radius $r>0$, define
\[
\mathcal{L}_{r,\theta}=\{x\in \nu L_{12}: \|x\|=r, \tau(x)=\theta\}.
\]
Letting $\mathcal{L}_{r,\theta}'$ be the corresponding object for $\fol'$, it is clear from the discussion above that $\mathcal{L}_{r,\theta}=\mathcal{L}_{r,\theta}'$ for any $\theta, r$. By equifocality of singular Riemannian foliations \cite{AT}, it follows that every leaf of $\fol$ is of the form $\exp(\mathcal{L}_{r,\theta})$ for some $r$ and $\theta$. Therefore, $\fol=\fol'$.
\end{proof}

\subsection{Concluding the proof}
We can finally gather all previous results in this section to provide the proof of the codimension 2 classification result.

\begin{proof}[Proof of Theorem~\ref{T:codim2class}]
	Let $(\CaP,\fol)$ be a codimension 2 polar foliation, with leaf space $X=\CaP/\fol$. By Proposition~\ref{P:triangles}, $X$ is isometric to a spherical triangle with angles $(\pi/2, \pi/2, \pi/g)$ for $g\in\{2,4, 6\}$ or $(\pi/2,\pi/3,\pi/4)$. In the first case, $g=6$ is excluded by the discussion at the beginning of \S\ref{subsec:22g}. Then the proof is a direct consequence of the previously stablished results:
	\begin{enumerate}[\rm(1)]
		\item If $X$ is a triangle with angles $(\pi/2, \pi/2, \pi/2)$, then $(\CaP,\fol)$ is the orbit foliation of the action of a Lie subgroup $\mathsf{Spin}_8\subset\mathsf{Spin}_9\subset\mathsf{F}_4$, by Proposition~\ref{prop:Spin8}.
		\item If $X$ is a triangle with angles $(\pi/2, \pi/2, \pi/4)$, then Theorem~\ref{T:pi/4} implies that $(\CaP,\fol)$ is either the orbit foliation of the action of $\Spin(7)\times \mathsf{SO}(2)\subset \Spin(9)\subset \mathsf{F}_4$ or of $\Spin(6)\times \Spin(3)\subset \Spin(9)\subset \mathsf{F}_4$, or congruent to the inhomogeneous foliation described in Proposition~\ref{P:nonhom-fol}.
		\item If $X$ is a triangle with angles $(\pi/2, \pi/3, \pi/4)$, then Theorem~\ref{T:irred} and Remark~\ref{rem:234} guarantee that $(\CaP,\fol)$ is the orbit foliation of the action of either $\mathsf{SU}(3)\cdot \mathsf{SU}(3)$ or $\mathsf{SO}(3)\cdot \mathsf{G}_2$.\qedhere
	\end{enumerate}	
\end{proof}

\bibliographystyle{alpha}
\bibliography{refs}

\end{document}